\documentclass[11pt,notitlepage,twoside,a4paper]{amsart}
 \usepackage{amsfonts}
\usepackage[latin1]{inputenc}
\usepackage[cyr]{aeguill}
\usepackage{xspace}
\usepackage[polutonikogreek,english]{babel}

\usepackage{amsmath,amssymb,enumerate}

\usepackage{epsfig,fancyhdr,color}

\usepackage{amssymb}
\usepackage{amsmath,amsthm}
\usepackage{latexsym}
\usepackage{amscd}
\usepackage{psfrag}
\usepackage{graphicx}
\usepackage[latin1]{inputenc}
\usepackage[all]{xy}
\usepackage[mathcal]{eucal}

\definecolor{NoteColor}{rgb}{1,0,0}

\renewcommand{\textsc}{\textcolor{red}}
\newcommand*{\tg}[1]{\textgreek{#1}}

\newtheorem*{theorem 1}{\rm\bf Proposition 1}
\newtheorem*{theorem 2}{\rm\bf Proposition 2}

\theoremstyle{definition}

\theoremstyle{remark}

\def\interieur#1{\mathord{\mathop{\kern 0pt #1}\limits^\circ}}

\title[Topology and biology]{Topology and biology: From Aristotle to  Thom}

\author{Athanase Papadopoulos}
\address{Athanase Papadopoulos,  Universit{\'e} de Strasbourg and CNRS,
7 rue Ren\'e Descartes,
 67084 Strasbourg Cedex, France, and Yau Mathematical Sciences Center, Tsinghua University, Jin Chun Yuan West Building, Haidian District, Beijing 100084, China}
\email{athanase.papadopoulos@math.unistra.fr}
\makeindex
 
\date{\today}

\begin{document}

 \begin{abstract}

   René Thom discovered several refined topological notions in the writings of Aristotle, especially the biological. More generally, he considered that some of the assertions of the Greek philosophers have a definite topological content, even if they were written 2400 years before the field of topology was born.
   In this article, we expand on these ideas of Thom. At the same time, we highlight the importance of biology in the works of Aristotle and Thom and we report on their conception of mathematics and more generally on science. 
   
   The final version of this article will appear in the book \emph{Geometry in history} (ed. S. G. Dani and A. Papadopoulos), Springer Verlag, 2019.
  \end{abstract}  
    \maketitle
  \medskip
   
   AMS classification: 03A05, 01-02, 01A20, 34-02, 34-03, 54-03, 00A30, 92B99.
   
   Keywords: Mathematics in Aristotle, topology in Aristotle, form, morphology, mathematics in biology, René Thom, morphogenesis,  hylemorphism, stratification, Stokes formula, history of topology.
   
  \tableofcontents

         \section{Introduction}

     Bourbaki,\index{Bourbaki} in his \emph{\'Eléments d'histoire des mathématiques}, at the beginning of the chapter on topological spaces, writes  (\cite[p.\,175]{Bourbaki}):
\begin{quote}\small
The notions of limit and continuity\index{continuous} date back to antiquity; one cannot make their complete history without studying systematically, from this point of view, not only the mathematicians, but also the Greek philosophers, and in particular Aristotle,\index{Aristotle (384--322 BCE)} and without tracing the evolution of these ideas through the mathematics of the Renaissance and the beginning of differential and integral calculus. Such a study, which would be interesting to carry on, would have to go much beyond the frame of the present note.\footnote{``Les notions de limite et de continuité remontent à l'antiquité ; on ne saurait en faire une histoire complète sans étudier systématiquement de ce point de vue, non seulement les mathématiciens, mais aussi les philosophes grecs et en particulier Aristote, ni non plus sans poursuivre l'évolution de ces idées à travers les mathématiques de la Renaissance et les débuts du calcul différentiel et intégral. Une telle étude, qu'il serait certes intéressant d'entreprendre, dépasserait de beaucoup le cadre de cette note."  [For all the translations that are mine, I have included the French text in footnote.]}
\end{quote}

The reference in this quote is to ``Greek philosophers" rather than ``Greek mathematicians." One may recall in this respect that at that time, mathematics was closely related to philosophy, even though there was a clear distinction between the two subjects. Many mathematicians were often philosophers and vice-versa, and mathematics was taught at the philosophical schools. A prominent example is Plato who, in his Academy, was essentially teaching mathematics. We may recall here a note by Aristoxenus,\index{Aristoxenus of Tarentum (b. c. 375 BCE)} the famous 4$^{\mathrm{th}}$ century BCE music theorist who, reporting on the teaching of Plato, writes that the latter used to entice his audience by announcing lectures on the Good or  similar philosophical topics, but that eventually the lectures turned out to be on geometry, number theory or astronomy, with the conclusion that ``the Good is One" \cite[p.\,98]{Aristoxene}. Aristotle\index{Aristotle (384--322 BCE)} was educated at Plato's academy, where he spent twenty years, learning mathematics in the pure Pythagorean tradition, and he was one of the best students there---Plato\index{Plato (428--348 BCE)} used to call him ``the Brain of the School"---. This is important to keep in mind, since Aristotle's mathematical side, which is the main topic of the present article, is usually neglected.

Over half a century since Bourbaki's statement, a systematic investigation of traces of topology\index{topology} in the works of the philosophers of  Greek antiquity has not been conducted yet. One problem for realizing such a program is that the ancient Greek texts that survive are generally analyzed and commented on by specialists in philosophy or ancient languages with little knowledge of mathematics, especially when it comes to topology, a field where detecting the important ideas requires a familiarity with this topic. In particular, it happens sometimes that the translation of terms conveying topological ideas do not render faithfully their content; this is the well-known problem of mathematical texts translated by non-mathematicians.
 
 Thom\index{Thom, René (1923--2002)} offers a rare, even unique, and important exception. He spent a great deal of his time reading (in Greek) the works of Aristotle\index{Aristotle (384--322 BCE)} and detecting the topological ideas they contain. He also corrected a number of mistakes in the existing translations and proposed more precise ones. In his 1991 article \emph{Matière, forme et catastrophes} \cite[p.\,367]{A-Matiere}, he declares that his past as a mathematician and initiator of catastrophe theory probably confers upon him a certain capacity of noticing in the work of Aristotle aspects that a traditional philosophical training would have put aside. He writes in his article \emph{Les intuitions topologiques primordiales de l'aristot\'elisme} \cite{Thom-Intuitions} (1988) that the earliest terms of Greek topology\index{topology} available in the Greek literature are contained in Parmenides'\index{Parmenides of Elea (6th-5th c. BCE)} poem, \emph{On Nature}. He notes, in particular, that the word \emph{sunekhes} (\tg{suneq'es}) which appears there is usually translated by \emph{continuous}, which is a mistake, because this word rather refers to (arcwise) connectedness. In order to make such a comment, and, especially given that no precise topological terminology existed in antiquity, a definite understanding of the topological ideas and intuitions involved together with their equivalent in modern language is required. There are many other examples of remarks of the same kind by Thom,\index{Thom, René (1923--2002)} and we shall note a few of them later in the text. Thanks to Thom's translation and interpretation of Aristotle\index{Aristotle (384--322 BCE)} and of other Greek philosophers, certain passages that were  obscure suddenly become clearer and more understandable for a mathematician. This is one of the themes that we develop in this article.

     In several of his works, Thom\index{Thom, René (1923--2002)} expressed his respect for the philosophers of Greek antiquity. He started to quote them extensively in his books and articles that he wrote in the beginning of the 1960s, a time when his mathematical work became part of a broader inquiry that included biology\index{biology} and philosophy. In the first book he published, \emph{Stabilité structurelle et morphogenèse :  Essai d'une th\'eorie g\'en\'erale des mod\`eles}, written in the period 1964--1968, he writes \cite[p.\,9]{SS}:\footnote{We quote from the English translation, \emph{Structural stability and morphogenesis. An outline of a general theory of models} published in 1975.}
     
\begin{quote}\small 
  If I have quoted the aphorisms of Heraclitus\index{Heraclitus (c. 544--480 BCE)} at the beginning of some chapters, the reason is that nothing else could be better adapted to this type of study. In fact, all the basic intuitive ideas of morphogenesis can be found in Heraclitus: all that I have done is to place these in a geometric and dynamic framework that will make them some day accessible to quantitative analysis. The ``solemn, unadorned words,"\footnote{Thom\index{Thom, René (1923--2002)} quotes here Heraclitus, Fragment No. 12: ``And the Sibyl with raving mouth, uttering words solemn, unadorned, and unsweetened, reaches with her voice a thousand years because of the god in her."  (Transl. A. Fairbanks, In: The First Philosophers of Greece, Scribner, 1898).} like those of the sibyl that have sounded without faltering throughout the centuries, deserve this distant echo.
\end{quote}

     Aristotle\index{Aristotle (384--322 BCE)}, to whom Thom referred constantly in his later works, placed biology\index{biology} and morphology\index{morphology} at the basis of almost every field of human knowledge, including physics, psychology, law and politics. Topology is involved, often implicitly but also explicitly, at several places in his writings.
  Bourbaki, in the passage cited at the beginning of this article, talks about limits and continuity\index{continuous}, but Thom\index{Thom, René (1923--2002)} highlighted many other fundamental topological questions in the works of the Stagirite. Among them is the central problem of topology,\index{topology} namely, that of reconstructing the global from the local. He noted that Aristotle\index{Aristotle (384--322 BCE)} saw that this problem, which he considered mainly from a morphological point of view, is solved in nature, especially in the vegetal world, where one can reconstruct a plant from a piece of it. Even in the animal world, where this phenomenon is exceptional---in general, one cannot resurrect an animal from one of its parts---he saw that if we amputate certain mollusks, crustaceans and insects of one of their organs, they can regenerate it. Most especially, Aristotle understood that it is possible to split certain eggs to obtain several individuals. 
               
       In the foreword to his   \emph{Esquisse d'une s\'emiophysique}, a book published in 1988, which carries the subtitle \emph{Physique aristot\'elicienne et th\'eorie des catastrophes},\footnote{We quote from the English translation published in 1990 under the title  
    \emph{Semio Physics: A Sketch}, with the subtitle  \emph{Aristotelian Physics and Catastrophe Theory.}} Thom\index{Thom, René (1923--2002)} writes \cite[p.\,{\sc viii}]{Esquisse}:\footnote{Page numbers refer to the English translation.} 
     ``It was only quite recently, almost by chance, that I discovered the work of Aristotle\index{Aristotle (384--322 BCE)}. It was fascinating reading, almost from the start." 
     More than a source of inspiration, Thom found in the writings of the Stagirite a confirmation of his own ideas.  At several occasions,  he expressed his surprise to see in these writings statements that he himself had came up with, in slightly different terms, without being aware of their existence. He declares in his 1991 article \emph{Matière, forme et catastrophes} \cite[p.\,367]{A-Matiere}: ``In the last years, in view of my writing of the \emph{Semiophysics}, I was led to go more thoroughly into the knowledge of the work of the Stagirite."\footnote{Ces dernières années, en vue de la rédaction de ma \emph{Sémiophysique}, j'ai dû entrer plus profondément dans la connaissance de l'\oe uvre du Stagirite.}   
  He found in that work the idea that, because we cannot understand everything in nature, we should concentrate our efforts on stable\index{stability} and generic\index{genericity} phenomena. Aristotle\index{Aristotle (384--322 BCE)} formulated the difference between a generic\index{genericity} and a non-generic phenomenon as a difference between a ``natural phenomenon" and an ``accident."\index{accident (Aristotle)} Interpreting this idea (and related ones) in terms of modern topology\index{topology} is one of the epistemological contributions of Thom.\index{Thom, René (1923--2002)} He writes in his \emph{Esquisse} \cite[p.\,12]{Esquisse} (p. ix of the English translation):  
 \begin{quote}\small
 [\ldots] If I add that I found in Aristotle\index{Aristotle (384--322 BCE)} the concept of genericity\index{genericity} (\tg{<ws >ep`i tò pol'u}), the idea of stratification\index{stratification!of an animal organism} as it might be glimpsed in the decomposition of the organism into homeomerous\index{homeomerous} and anhomeomerous\index{anhomeomerous} parts by Aristotle the biologist, and the idea of the breaking down of the genus into species as images of bifurcation, it will be agreed that there was matter for some astonishment.
 \end{quote}
 
 The notion of ``homeomerous" (\tg{<omoiomer'hs}), to which Thom refers in this passage, occurs in Aristotle's\index{Aristotle (384--322 BCE)} \emph{History of animals}, his \emph{Parts of animals}, and in other zoological writings. It contains the idea of self-similarity that  appears in mathematics. Among the other topological ideas that Thom\index{Thom, René (1923--2002)} found in Aristotle's\index{Aristotle (384--322 BCE)} writings are those of open and closed set, of cobordism,\index{cobordism} and the Stokes formula.\index{Stokes formula}  We shall quote the explicit passages  below. This Aristotelian view of mathematics, shared by Thom,\index{Thom, René (1923--2002)} tells us in particular that mathematics emerges from our daily concepts, rather than exists is some Platonist ideal realm.

 Thom published several articles on topology\index{topology} in Aristotle\index{Aristotle (384--322 BCE)}, including \emph{Les intuitions topologiques primordiales de l'aristot\'elisme}  (\emph{The primary intuitions of Aristotelianism}) \cite{Thom-Intuitions} (1988) which is an expanded version of a section in Chapter 7 of the \emph{Esquisse}, \emph{Matière, forme et catastrophes} (\emph{Matter, form and catastrophes})  \cite{A-Matiere} (1991), an article published in the proceedings of a conference held at the Unesco headquarters in Paris at the occasion of the 23$^{\mathrm{rd}}$ centenary of the Philosopher, \emph{Aristote topologue} (Aristotle as a topologist) \cite{Aristote-Topologue} (1999) and others. It comes out from Thom's articles that even though Aristotle did not build any topological theory in the purely mathematical sense, topological intuitions are found throughout his works. Thom\index{Thom, René (1923--2002)} writes in  \cite[p.\,395]{Thom-Intuitions}: ``We shall not find these intuitions in explicit theses and constructions of the theory. We shall find them mostly in some `small sentences' that illuminate the whole corpus with their brilliant concision."\footnote{Ces intuitions, on ne les trouvera pas dans les thèses et les constructions explicites de la théorie. On les trouvera surtout dans quelques ``petites phrases", qui illuminent tout le corpus de leur éclatante concision.}

  Thom\index{Thom, René (1923--2002)} considered that mathematics is a language that is adequate for philosophy. he expressed this is several papers, and in particular in the article \emph{Logos phenix} \cite{T-Logos} and others that are reprinted in the collection \emph{Modèles mathématiques de la morphogenèse} \cite{Morpho}.  In fact, beyond the discussion that we have here concerning topological notions, the question of how much philosophical notions carrying a mathematical name (infinity, limit, continuity, etc.) differ from their mathematical counterparts is debatable. Indeed, one may argue that mathematics needs precise terms and formal definitions. Such a position is expressed by Plotnitsky, who writes in his essay published in the present volume \cite{Plotnitsky}: 
  \begin{quote}\small
  What made topology a mathematical discipline is that one can associate algebraic structures (initially numbers, eventually groups and other abstract algebraic structures, such as rings) to the architecture of spatial objects that are invariant under continuous transformations, independently of their geometrical properties.
   \end{quote}
Plotnitsky also quotes Hermann Weyl\index{Weyl, Hermann (1885--1955)} from his book \emph{The continuum} \cite{Weyl}, where the latter says that the concepts offered to us by the intuitive notion of continuum cannot be identified with those that mathematics presents to us. Thom had a different point of view. For him, mathematics and our intuition of the real world are strongly intermingled. In the article \emph{Logos phenix} \cite[p.\,292]{Morpho}, he writes: ``How can we explain that mathematics represents the real world? The answer, I think is offered to us by the intuition of the continuous\index{continuous}"\footnote{Comment expliquer que les mathématiques représentent le réel ? La réponse, je crois, nous est offerte par l'intuition du continu."} For him, it is the geometrical notion of continuous\index{continuous} that gives a meaning to beings that would have needed infinitely many actions. He mentions the paradox of Achilles and the Tortoise that allows us to give a meaning to the infinite sum 1/2+1/4+1/8+... Here, he writes, ``infinite becomes seizable in action."\footnote{L'infini devient saisissable en acte.} Conversely, he says, the introduction of a continuous\index{continuous} underlying substrate allows us to explain the significant---non-trivial---character of many mathematical theorems  \cite[p.\,293]{Morpho}.
  
  Beyond philosophy, Thom\index{Thom, René (1923--2002)} used the mathematics that he developed in order to express natural phenomena.  At several occasions, he insisted on the fact that the mathematical concepts need not be formalized or rigorously expressed in order to exist. He recalls in \cite{Atiyah} that ``rigor" is a latin word that reminds him of the sentence \emph{rigor mortis} (the rigor of a dead body). He writes that rigor is ``a very unnecessary quality in mathematical thinking." In his article \cite{Thom-enc}, he declares that anything which is rigorous is not significant.\footnote{Tout ce qui est rigoureux est insignifiant.}     
  
  Our aim in the next few pages is to expand on these ideas.  A collection of essays on the Thom and his work appeared in print recently, \cite{Thom}.

 \section{D'Arcy Thompson}
 Before talking more thoroughly about Aristotle\index{Aristotle (384--322 BCE)} and Thom, I would like to say a few words on the mathematician, biologist and philosopher who, in many ways, stands between them, namely, D'Arcy Thompson.\index{Thompson, D'Arcy Wentworth (1860-1948)!\emph{On growth and form}}\footnote{D'Arcy Wentworth Thompson\index{Thompson, D'Arcy Wentworth (1860-1948)} (1860-1948) was a Scottish biologist with a profound passion for mathematics and for Greek science and philosophy.  He was also an accomplished writer and his \emph{magnum opus}, \emph{On growth and form}, is an authentic literary piece. He is considered as the first who found a relation between the Fibonacci sequence and some logarithmic spiral structures in the animal and vegetable life (mollusk shells, ruminant horns, etc.). His works were influential on several 20$^{\mathrm{th}}$ century thinkers including C. H. Waddington, Alan Turing, Claude Lévi-Strauss and Le Corbusier and on artists such as Richard Hamilton and Eduardo Paolozzi. D'Arcy Thompson\index{Thompson, D'Arcy Wentworth (1860-1948)} translated Aristotle's \emph{History of animals}, a translation which is among the authoritative ones.} At several occasions, Thom\index{Thom, René (1923--2002)} referred to the latter's book \emph{On growth and form}\index{Thompson, D'Arcy Wentworth (1860-1948)!\emph{On growth and form}}\index{Thompson, D'Arcy Wentworth (1860-1948)!\emph{On growth and form}} (1917), a frequently quoted work whose main object is the existence of mathematical models for growth and form in animal and vegetable biology\index{biology} at the level of cells, tissues and other parts of a living organism. D'Arcy Thompson\index{Thompson, D'Arcy Wentworth (1860-1948)} emphasized the striking analogies between the mathematical patterns that describe all these parts, searching for general laws based on these patterns. In reading \emph{On growth and form}\index{Thompson, D'Arcy Wentworth (1860-1948)!\emph{On growth and form}}, Thom\index{Thom, René (1923--2002)} had the same reaction he had when he read Aristotle\index{Aristotle (384--322 BCE)}: he was amazed by the richness of the ideas contained in it, and by their closeness to his own ideas. In turn, Thompson\index{Thompson, D'Arcy Wentworth (1860-1948)} had a boundless admiration for Aristotle who had placed biology at the center of his investigations.  In an article titled \emph{On Aristotle as biologist}, Thompson\index{Thompson, D'Arcy Wentworth (1860-1948)} calls the latter ``the great biologist of Antiquity, who is \emph{maestro di color che sanno},\footnote{The quote is from Dante's \emph{Inferno}, in which he mentions Aristotle\index{Aristotle (384--322 BCE)}: ``I saw the master of those who know."} in the science as in so many other departments of knowledge." \cite[p.\,11]{Thompson2} Likewise, Thom, in his articles \emph{Les intuitions topologiques primordiales de l'aristot\'elisme} and \emph{Matière, forme et catastrophes}, refers to Aristotle as ``the Master" (``le Maître"). Thompson writes in \emph{On growth and form}\index{Thompson, D'Arcy Wentworth (1860-1948)!\emph{On growth and form}} \cite[p.\,15]{Thompson2}: 
\begin{quote}\small
[Aristotle] recognized great problems of biology\index{biology} that are still ours today, problems of heredity, of sex, of nutrition and growth, of adaptation, of the struggle for existence, of the orderly sequence of Nature's plan. Above all he was a student of Life itself. If he was a learned anatomist, a great student of the dead, still more was he a lover of the living. Evermore his world is in movement. The seed is growing, the heart is beating, the frame breathing. The ways and habits of living things must be known: how they work and play, love and hate, feed and procreate, rear and tend their young; whether they dwell solitary, or in more and more organized companies and societies. All such things appeal to his imagination and his diligence.
\end{quote}

Before Thom, D'Arcy Thompson\index{Thompson, D'Arcy Wentworth (1860-1948)} foresaw the importance of topology\index{topology} in biology.\index{biology} He writes, in the same treatise \cite[p.\,609 ff]{Thompson2}:
  \begin{quote}\small
  [\ldots] for in this study of a segmenting egg we are on the verge of a subject adumbrated by Leibniz, studied more deeply by Euler, and greatly developed of recent years. [\ldots]  Topological analysis seems somewhat superfluous here; but it may come into use some day to describe and classify such complicated, and diagnostic, patterns as are seen in the wings of a butterfly or a fly. 
  \end{quote}

 A question which is dear to topologists and which we already mentioned is that of the ``local implies global."\index{local and global} Thompson\index{Thompson, D'Arcy Wentworth (1860-1948)} addressed it in his book
\emph{On growth and form}\index{Thompson, D'Arcy Wentworth (1860-1948)!\emph{On growth and form}}. He writes (p. 2019):
\begin{quote}\small
The biologist, as well as the philosopher, learns to recognize that the whole is not merely the sum of its parts. It is this, and much more than this. For it is not a bundle of parts but an organization of parts, of parts in their mutual arrangement, fitting one with another, in what Aristotle\index{Aristotle (384--322 BCE)} calls ``a single and indivisible principle of unity"; and this is no merely metaphysical conception, but is in biology\index{biology} the fundamental truth which lies at the basis of Geoffroy's (or Goethe's) law of ``compensation," or ``balancement of growth."
\end{quote}
 
It is of course not coincidental that Thompson\index{Thompson, D'Arcy Wentworth (1860-1948)} mentions Goethe\index{Goethe, Johann Wolfgang von (1749--1832)} and  Geoffroy-Saint-Hilaire \index{Geoffroy-Saint-Hilaire, Isidore (1805--1861)} and it is well worth to recall here that Goethe\index{Goethe, Johann Wolfgang von (1749--1832)} was also a passionate student of biology, in particular, of form\index{form} and morphology.\index{morphology} He wrote an essay on the evolution of plants based on their form. To him is generally attributed the first use of the word ``metamorphosis"\index{metamorphosis} in botanics---even though the concept was present in Aristotle, who studied the metamorphosis of butterflies, gnats and other insects---. 
             
  Isidore Geoffroy Saint-Hilaire,\index{Geoffroy-Saint-Hilaire, Isidore (1805--1861)} the other  naturalist to whom D'Arcy Thompson refers, wrote an important treatise in three volumes titled \emph{Histoire naturelle générale des règnes organiques} (General natural history of organic worlds) in which he carries out a classification of animals that is still used today. In commenting on Aristotle\index{Aristotle (384--322 BCE)}, Goeffroy writes in this treatise \cite[vol. I, p. 19 ss.]{Geoffroy}:
  \begin{quote}\small
   He is, in every branch of human knowledge, like a master who develops it alone. He reaches, he extends, the limit of all sciences, and at the same time he penetrates their intimate depths. From this point of view, Aristotle is an absolutely unique exception in the history of human thought, and if something here is amazing, it is not the fact that this exception is unique, but that there exists one, as such a meeting of faculties and of knowledge is surprising for anyone who wants to notice it psychologically. [\ldots] Among his multiple treatises, the \emph{History of animals} and the \emph{Parts} are the main monuments of his genius.\footnote{Il est, dans chaque branche du savoir humain, comme un maître qui la cultiverait seule; il atteint, il recule les limites de toutes les sciences, et il en pénètre en même temps les profondeurs intimes. Aristote est, à ce point de vue, une exception absolument unique dans l'histoire de l'esprit humain, et si quelque chose doit nous étonner ici, ce n'est pas qu'elle soit restée unique, c'est qu'il en existe une : tant une semblable réunion de facultés et de connaissances est surprenante pour qui veut s'en rendre compte psychologiquement. [\ldots] Entre ses nombreux traités, les deux monuments principaux de son génie sont l'\emph{Histoire des animaux} et \emph{le Traité des parties}.]}
   \end{quote}

To close this section, let me mention that the biologist Thomas Lecuit, newly appointed professor at Collège de France, and gave his first course, in the year 2017-2018, on the problem of morphogenesis, started his first lecture by a tribute to D'Arcy Thompson, emphasizing his contribution to the problem of understanding the diversity of forms\index{form} and the mathematical patterns underlying them and their transformations, and highlighting the continuing importance of his 100 years old book  \emph{On growth and form}\index{Thompson, D'Arcy Wentworth (1860-1948)!\emph{On growth and form}}. Let me also mention that in 2017, a workshop dedicated to that book was held at Princeton's Institute for Advanced Study.

   \section{Aristotle\index{Aristotle (384--322 BCE)}, mathematician and topologist}

   Aristotle\index{Aristotle (384--322 BCE)} is very poorly known as a mathematician, a situation which is unfair. Although he did not write any mathematical treatise (or, rather, there is no indication that such a treatise existed), Aristotle had an enormous influence on mathematics, by his thorough treatment of first principles, axioms, postulates and other foundational notions of geometry, by his classification of the various kinds of logical reasonings, his deep thoughts on the use of motion (that is, what we call today isometries) in axioms and in proofs, and on the consequences of that use. There are may ways in which Aristotle preceded Euclid, not in compiling a list of axioms, but in his profound vision on the axiomatic approach to geometry. We refer in particular to his discussions of first principles in his \emph{Posterior analytics}\index{Aristotle (384--322 BCE)!\emph{Posterior analytics}}\footnote{\emph{Posterior analytics}, \cite{A-Posterior} 74b5 and 76a31-77a4.}, to his insistence in the \emph{Metaphysics}\index{Aristotle (384--322 BCE)!\emph{Metaphysics}}\footnote{\emph{Metaphysics} \cite{A-Metaphysics} 997a10.} on the fact that a demonstrative science is based on axioms that are not provable, to his discussion of the \emph{reductio ad absurdum} reasoning in the same work,\footnote{\emph{Posterior analytics}, \cite{A-Posterior}  85a16ff.} and there are many other ideas on the foundations of mathematics in his work that we could have mentioned. Besides that, several mathematical propositions are found throughout his works\index{Aristotle (384--322 BCE)}. We find in particular results related to all the fundamental problems of mathematics of that epoch: the parallel problem in the \emph{Prior analytics}\index{Aristotle (384--322 BCE)!\emph{Prior analytics}},\footnote{\emph{Prior analytics} \cite{A-Prior} 65a4--9, 66a11-15.} the incommesurability\index{incommensurability} of the diagonal of a square\footnote{\emph{Prior analytics} \cite{A-Prior} 65b16-21.} and the squaring of the circle\index{squaring of the circle} (by means of the squaring of lunules) in the same treatise,\footnote{\emph{Prior analytics} \cite{A-Prior} 69a20-5.}  etc.  There is an axiom in the foundations of geometry that was given the name \emph{Aristotle's axiom};\index{Aristotle's axiom} see Greenberg's article \cite{Greenberg}. In Book III of 
   \emph{On the heavens},\index{Aristotle (384--322 BCE)!\emph{On the heavens}} talking about form,\index{form} one of his favorite topics, and commenting on a passage of Plato's\index{Plato (428--348 BCE)} \emph{Timaeus}\index{Plato (428--348 BCE)!\emph{Timaeus}} in which regular polyhedra are associated with the four sublunar elements (earth, water, air, fire), Aristotle\index{Aristotle (384--322 BCE)} states that there are exactly three regular figures that tile the plane, namely, the equilateral triangle, the square and the regular hexagon, and that in space, there are only two: the pyramid and the cube.\footnote{\emph{On the heavens} \cite{A-Heavens} 306b1-5.} It is possible that this passage is the oldest surviving written document in which this theorem is stated.\footnote{The reason why Aristotle makes this statement here is not completely clear, but it is reasonable to assume that it is because Plato\index{Plato (428--348 BCE)}, in the \emph{Timaeus}---a  text which is very mysterious---conjectured that the elementary particles of the four elements have the form of the regular polyhedra he associated  to them; hence the question raised by Aristotle\index{Aristotle (384--322 BCE)} concerning the  tiling of space using regular polyhedra.}
        One may also consider Aristotle's \emph{Problems},\index{Aristotle (384--322 BCE)!\emph{Problems}}\footnote{The translations below are from \cite{A-Problems}.} a treatise in 38 books, assembled by themes, one of the longest of the Aristotelian corpus. It consists of a list of commented (open) problems, of the kind mathematicians are used to edit, with the difference that Aristotle's problems concern not only mathematics but all subjects of human knowledge: psychology, biology, physics, acoustics, medicine, ethics, physiology, etc. In total, there are 890 problems. They generally start with the words ``Why is it that..."  For instance, Problem 48 of Book XXVI asks: \emph{Why is it that the winds are cold, although they are due to movement caused by heat?}  Problem 19 of Book XXXI asks: \emph{Why is it that when we keep our gaze fixed on objects of other colours our vision deteriorates, whereas it improves if we gaze intently on yellow and green objects, such as herbs and the like?} Problem 20 of the same book asks: \emph{Why is it that we see other things better with both eyes, but we can judge of the straightness of lines of writing better with one eye?} Some of the problems concern mathematics. For instance, Problem 3 of Book XV asks: \emph{Why do all men, Barbarians and Greek alike, count up to 10 and not up to any other number, saying for example, 2, 3, 4, 5 and then repeating then, ``one-five", ``two-five", just as they say eleven, twelve?}
         Problem 5 of the same book starts with the question: \emph{Why is it that, although the sun moves with uniform motion, yet the increase and decrease of the shadows is not the same in any equal period of time?}         Book XVI is dedicated to ``inanimate things." Problems 1 and 2 in that book concern floating bubbles. Problem 2 asks: \emph{Why are bubbles hemispherical?}\index{Thompson, D'Arcy Wentworth (1860-1948)!\emph{On growth and form}}\footnote{D'Arcy Thompson\index{Thompson, D'Arcy Wentworth (1860-1948)} was also fascinated by the questions of form\index{form} and transformation of floating bubbles,\index{floating bubble} in relation with the question of growth of a living cell submitted to a fluid pressure (cf. \emph{On growth and form}, Chapters V to VII). On p. 351, he talks about ``the peculiar beauty of a soap-bubble, solitary or in collocation [..] The resulting form\index{form} is in such a case so pure and simple that we come to look on it as well-nigh a mathematical abstraction." On p. 468, he writes: ``Bubbles have many beautiful properties besides the more obvious ones. For instance, a floating bubble is always part of a sphere, but never more than a hemisphere; in fact it is always rather less, and a very small bubble is considerably less than a hemisphere. Again, as we blow up a bubble, its thickness varies inversely as the square of its diameter; the bubble becomes a hundred and fifty times thinner as it grows from an inch in a diameter to a foot."  Later in the same chapter, Thompson talks about clustered bubbles (p. 485). He quotes Plateau on soap-bubble shapes.} Problem 5 asks: \emph{Why is it that a cylinder, when it is set in motion, travels straight and describes straight lines with the circles in which it terminates, whereas a cone revolves in a circle, its apex remaining still, and describes a circle with the circle in which it terminates?}  Problem 5 of the same book concerns the traces of oblique sections of a cylinder rolling on a plane.  Problem 6 concerns a property of straight lines: \emph{Why is it that the section of a rolled book, which is flat, if you cut it parallel to the base becomes straight when unrolled, but if it is cut obliquely becomes crooked?}

Another mathematical topic discussed in detail in Aristotle's\index{Aristotle (384--322 BCE)} works is that of infinity,\index{infinity} for which the Greeks had a name, \emph{apeiron} (\tg{>'apeiron}), meaning ``boundless." This notion, together with that of limit, is discussed in the \emph{Categories},\index{Aristotle (384--322 BCE)!\emph{Categories}} the \emph{Physics},\index{Aristotle (384--322 BCE)!\emph{Physics}} the \emph{Metaphysics}\index{Aristotle (384--322 BCE)!\emph{Metaphysics}}, \emph{On the heavens},\index{Aristotle (384--322 BCE)!\emph{On the heavens}} and in other writings.  In the \emph{Physics},\footnote{\emph{Physics}, \cite{A-Physics} 201a-b.} Aristotle\index{Aristotle (384--322 BCE)} mentions the two occurrences of the infinite in mathematics: the infinitely large, where, he says, ``every magnitude is surpassed"  and the infinitely small, where ``every assigned magnitude is surpassed in the direction of smallness."  Aristotle\index{Aristotle (384--322 BCE)} disliked the ``unbounded infinite" and his interest lied in the second kind. He writes in the same passage:  
  \begin{quote}\small
  Our account does not rob the mathematicians of their science, by disproving the actual existence of the infinite in the direction of increase, in the sense of the untraversable. In point of fact they do not need the infinite and do not use it. They postulate only that the finite straight line may be produced as far as they wish.
  \end{quote}
  This is a reference to the occurence of the infinitely large in the axioms of geometry (predating Euclid). It is interesting to note that Descartes\index{Descartes, René (1596--1650)} considered that only the infinitely small appears in mathematics. The infinitely large, in his conception,  belongs to metaphysics only. See the comments by R. Rashed in his article \emph{Descartes et l'infiniment petit} \cite{Rashed-Descartes}.

 Aristotle\index{Aristotle (384--322 BCE)}, in his discussion of infinity,\index{infinity} makes a clear distinction between the cases where infinity is attained or not.  In a passage of the \emph{Physics}\index{Aristotle (384--322 BCE)!\emph{Physics}},\footnote{\emph{Physics}, \cite{A-Physics} 204a.}  he provides a list of the various senses in which the word ``infinite" is used: (1) infinity incapable of being gone through; (2)   infinity capable of being gone through having no termination; (3) infinity  that ``scarcely admits of being gone through";  (4) infinity  that ``naturally admits of being gone through, but is not actually gone through or does not actually reach an end."  He discusses the possibility for an infinite body to be simple infinite or compound infinite.\footnote{\emph{Physics}, \cite{A-Physics} 204b10.} In the same page, he talks about form,\index{form} which, he says, ``contains matter and the infinite."  There is also a mention of infinite series in the \emph{Physics}.\footnote{\emph{Physics} \cite{A-Physics} 206a25-206b13.}

          One should also talk about the mathematical notion of continuity in the writings of the Philosopher.
          
            In the \emph{Categories}, Aristotle\index{Aristotle (384--322 BCE)} starts by classifying  quantities into discrete\index{discrete} or continuous.\index{continuous!vs. discrete}\index{discrete!vs. continuous}\index{continuous}\index{Aristotle (384--322 BCE)!\emph{Categories}}\footnote{\emph{Categories} \cite{A-Categories} 4b20.} He declares that some quantities are such that ``each part of the whole has a relative position to the other parts; others have within them no such relation of part to part," a reference to topology.\index{topology} As examples of discrete\index{discrete} quantities, he mentions number (the integers) and speech.  As examples of continuous\index{continuous} quantities, he gives lines, surfaces, solids, time and place\index{place} (a further reference to topology). He explains at length, using the mathematical language at his disposal, why the set of integers is discrete.\index{discrete} According to him, two arbitrary integers ``have no common  boundary,\index{boundary} but are always separate." He declares that the same holds for speech: ``there is no common boundary at which the syllables join, but each is separate and distinct from the rest." This is a way of saying that each integer (respectively each syllable) is isolated from the others. Aristotle\index{Aristotle (384--322 BCE)} writes that a line is a continuous\index{continuous} quantity  ``for it is possible to find a common boundary at which its parts join."   He decalres that space is a continuous quantity ``because the parts of a solid occupy a certain space, and these have a common boundary; it follows that the parts of space also, which are occupied by the parts of the solid, have the same common boundary as the parts of the solid." Time, he writes, is also a continuous quantity, ``for its parts have a common boundary." The notion of boundary\index{boundary} is omnipresent in this discussion. Thom\index{Thom, René (1923--2002)} was fascinated by this fact and he emphasized it in writings. We shall discuss this in more detail later in this paper.

In another passage of the \emph{Categories}, Aristotle\index{Aristotle (384--322 BCE)} discusses position and relative position, notions that apply to both discrete\index{discrete} and continuous\index{continuous} quantities:  ``Quantities consist either of parts which bear a relative position each to each, or of parts which do not."\index{Aristotle (384--322 BCE)!\emph{Categories}}\footnote{Categories \cite{A-Categories} 5a10.} Among the quantities of the former kind, he mentions the line, the plane,  the solid space, for which one may state ``what is the position of each part and what sort of parts are contiguous." On the contrary, he says, the parts of the integers do not have any relative position each to each, or a particular position, and it is impossible to state what parts of them are contiguous. The parts of time, even though the latter is a continuous\index{continuous} quantity, do not have position, because he says, ``none of them has an abiding existence, and that which does not abide can hardly have position." Rather, he says, such parts have a \emph{relative order}, like for number, and the same holds for speech: ``None of its parts has an abiding existence: when once a syllable is pronounced, it is not possible to retain it, so that, naturally, as the parts do not abide, they cannot have position." In this and in other passages of Aristotle's\index{Aristotle (384--322 BCE)} work,
motion and the passage of time are intermingled with spatiality. It is interesting to see that Hermann Weyl,\index{Weyl, Hermann (1885--1955)} in his book \emph{Space, time and matter},  also insisted on the importance of the relation between, on the one hand, motion, and, on the other hand, space, time and matter  \cite[p.\,1]{Weyl-STM}: ``It is the composite idea of \emph{motion}  that these three fundamental conceptions enter into intimate relationship."

  Talking about position, we come to the important notion of place.
  
Several Greek mathematicians insisted on the difference between space and place, and Aristotle\index{Aristotle (384--322 BCE)} was their main representative. They used the word khôra  (\tg{q'wra})\index{topos and chora} for the former and topos (\tg{t'opos}) for the latter. Again, Thom regards Aristotle's\index{Aristotle (384--322 BCE)} discussion from a mathematician point of view, and he considers it as readily leading to a topological mathematization,
 although it does not use the notational apparatus (which, at that time, was nonexistent) or the technical language of topology\index{topology} to which we are used today. In the \emph{Physics}\index{Aristotle (384--322 BCE)!\emph{Physics}}, Aristotle\index{Aristotle (384--322 BCE)} gives the following characteristics of place\footnote{\emph{Physics}, \cite{A-Physics} 211a ff.}: (1) Place is what contains that of which it is the place.\index{place}
(2) Place is no part of the thing.
(3) The immediate place of a thing is neither less nor greater than the thing.
(4) Place can be left behind by the thing and is separable. (5) All place admits of the distinction of up and down, and each of the bodies is naturally carried to its appropriate place and rests there, and this makes the place either up or down.
   
This  makes Aristotle's\index{Aristotle (384--322 BCE)} concept of place\index{place} close to our mathematical notion of boundary.\index{boundary} Furthermore, it confers to the notion of place the status of a ``relative" notion: a place is defined in terms of boundary, and the boundary is also the boundary of something else. At the same time, this is not too far from the notion of ``relative position" that was formalized later on in Galilean mechanics.\index{Galilean physics}
Aristotle writes in the same passage:
\begin{quote}\small
 We ought to try to make our investigations such as will render an account of place, and will not only solve the difficulties connected with it, but will also show that the attributes supposed to belong to it do really belong to it, and further will make clear the cause of the trouble and of the difficulties about it.
\end{quote}

Aristotle then introduces a dynamical aspect in his analysis of place. He discusses motion and states that locomotion and the phenomena of increase and diminution involve a variation of place. In the same passage of the \emph{Physics}\index{Aristotle (384--322 BCE)!\emph{Physics}}, he writes:
  \begin{quote}\small
  When what surrounds, then, is not separate from the thing, but is in continuity with it, the thing is said to be in what surrounds it, not in the sense of in place, but as a part in a whole. But when the thing is separate or in contact, it is immediately `in' the inner surface of the surrounding body, and this surface is neither a part of what is in it nor yet greater than its extension, but equal to it; for the extremities of things which touch are coincident.
  \end{quote}

Regarding place\index{place} and its relation to boundary,\index{boundary} we mention another passage from the \emph{Physics}\index{Aristotle (384--322 BCE)!\emph{Physics}}:\footnote{\emph{Physics}, \cite{A-Physics} 211b5-9.} ``Place is [\ldots] the boundary of the containing body at which it is in contact with the contained body. (By the contained body is meant what can be moved by way of locomotion)." Thom\index{Thom, René (1923--2002)} commented on this passage at several occasions. In his paper \emph{Aristote topologue} \cite{Aristote-Topologue}, he discusses the relation, in Aristotle's writings, between topos\index{topos and eschata} and \emph{eschata} (\tg{>'esqata}), that is, the limits, or extreme boundaries.  We shall elaborate on this below.

The overall discussion of place in Aristotle's\index{Aristotle (384--322 BCE)}\index{place} work, and its relation with shape and boundary\index{boundary} is involved, and the various translations of the relevant passages in his writings often differ substantially from each and depend on the understanding of their author.  This not a surprise, and the reader may imagine the difficulties in formulating a precise definition of boundary which does not use the language of modern topology.\index{topology} As a matter of fact, Aristotle emphasized the difficulty of defining place.  In the \emph{Physics}\index{Aristotle (384--322 BCE)!\emph{Physics}}, he  writes:  
\begin{quote}\small Place\index{place} is thought to be something important and hard to grasp, both because the matter and the shape present themselves along with it, and because the displacement of the body that is moved takes place in a stationary container, for it seems possible that there should be an interval which is other than the bodies which are moved.   [\dots] 
 Hence we conclude that \emph{the innermost motionless boundary\index{boundary} of what contains is place.} [\dots] Place\index{place} is thought to be a kind of surface, and as it were a vessel, i.e. a container of the thing. Further; place is coincident with the thing, for boundaries are coincident with the bounded.\footnote{\emph{Physics} \cite{A-Physics} 212a20.}
    \end{quote}

 Place\index{place} is also strongly related to form.\index{form}  Aristotle\index{Aristotle (384--322 BCE)} states that form is the boundary\index{boundary} of the thing whereas place\index{place} is the boundary of the body that contains  the thing. 
Thom interpreted the texts where Aristotle\index{Aristotle (384--322 BCE)} makes a distinction between ``the boundary of a thing" and ``the boundary of the body that contains it"---a structure of ``double ring" of the eschata,\index{topos and eschata} as he describes it---using homological considerations. More precisely, he saw there a version of the Stokes formula.\index{Stokes formula} We shall review this in \S \ref{s:Thom} below.

 Let us quote a well-known text which belongs to the Pythagorean literature, written by Eudemus of Rhodes (4th c. BCE), in which the latter quotes Archytas of Tarentum\index{Archytas of Tarentum (428--347 BCE)}. This text shows the kind of questions on space and place that the Greek philolophers before Aristotle addressed, e.g., whether space is bounded or not, whether ``outer space" exists, and the paradoxes to which this existence leads (see \cite[p.\,541]{Huffman}):
\begin{quote}\small
``But Archytas," as Eudemus says, ``used to propound the argument
in this way: `If I arrived at the outermost edge of the
heaven [that is to say at the fixed heaven], could I extend my
hand or staff into what is outside or not?' It would be paradoxical
not to be able to extend it. But if I extend it, what is outside
will be either body or place. It doesn't matter which, as we will
learn. So then he will always go forward in the same fashion to
the limit that is supposed in each case and will ask the same
question, and if there will always be something else to which
his staff [extends], it is clear that it is also unlimited. And if it is
a body, what was proposed has been demonstrated. If it is
place, place\index{place} is that in which body is or could be, but what is
potential must be regarded as really existing in the case of eternal things, and thus there would be unlimited body and space."
(Eudemus, Fr. 65 Wehrli, Simplicius, In Ar. Phys. iii 4; 541)
\end{quote}

Let \index{Simplicius, VIth c. CE} us now pass to other aspects of mathematics in Aristotle\index{Aristotle (384--322 BCE)}.   We mentioned in the introduction  his notion of homeomerous.\index{homeomerous} This is used at several occasions in his works, and especially in his zoological treatises. He introduces it at the beginning of the \emph{History of animals} \cite{A-HA}, where he talks about simple\index{simple part (Aristotle)} and complex parts.\index{complex part (Aristotle)} A part is simple, he says, if, when divided, one recovers parts that have the same form as the original part, otherwise, it is complex. For instance, a face is subdivided into eyes, a nose, a mouth, cheeks, etc., but not into faces. Thus, a face is a complex part of the body. On the other hand, blood, bone, nerves, flesh, etc. are simple parts because subdividing them gives blood, bone, nerves, flesh, etc.  Anhomeomorous\index{anhomeomerous} parts in turn are composed of homeomerous\index{homeomerous} parts: for instance, a hand is constituted of flesh, nerves and bones. The classification goes on. Among the  parts, some are called ``members." These are the parts which form a complete whole but which also contain distinct parts: for example, a head, a leg, a hand, a chest, etc. Aristotle also makes a distinction between parts responsible of ``act", which in general are inhomogeneous (like the hand) and the others, which are homogeneous (like blood) and which he considers as  ``potential parts."  
In D'Arcy Thompson's translation, the word ``homeomerous" is translated  by ``uniform with itself," and sometimes by ``homogeneous."  We mention that the term homeomoerous was also used in geometry, to denote 
lines that are self-similar in the sense that any part of them can be moved to coincide with any other part. Proclus,\index{Proclus of Lycia (411--485)} the 5$^{\mathrm{th}}$ century mathematician, philosopher and historian of mathematics, discusses this notion in his \emph{Commentary on the First Book of Euclid's Elements}.\index{Proclus of Lycia (411--485)!\emph{In Euclidem}} He gives as an example of a homeomerous\index{homeomerous} curve\index{homeomerous!curve}\index{curve!homeomerous} the cylindical helix,\index{cylindrical helix}\index{helix!cylindrical} attributing the definition to Apollonius\index{Apollonius of Perga (c. 262--c. 190 BCE)} in his book titled \emph{On the screw}\index{Apollonius of Perga (c. 262--c. 190 BCE)!\emph{On the screw}} and which does not survive. Among the curves of similar shape that are not homeomerous,\index{anhomeomerous}\index{anhomeomerous!curve}\index{curve!anhomeomerous} he gives the example of Archimedes' (planar) spiral,\index{Archimedes spiral} the conical helix,\index{conical helix}\index{helix!conical} and the spherical helix\index{spherical helix}\index{helix!spherical} (see p. 95 of ver Eecke's edition of his \emph{Commentary on the First Book of Euclid's Elements} \cite{Proclus}).\index{Proclus of Lycia (411--485)!\emph{In Euclidem}}  In the same work, Proclus\index{Proclus of Lycia (411--485)} attributes to Geminus\index{Geminus of Rhodes (1$^\mathrm{st}$ c. BCE)} a result stating that there are only three homeomerous curves:\index{homeomerous!curve} the straight line, the circle and the cylindrical helix \cite[p.\,102]{Proclus}.

 Aristotle\index{Aristotle (384--322 BCE)} was a mathematician in his unrelenting desire for making exhaustive classifications and of finding \emph{structure} in phenomena: after all, he introduced his \emph{Categories}  for that purpose.   In the \emph{History of animals}, like in several other treatises by Aristotle\index{Aristotle (384--322 BCE)}, the sense of detail is dizzying. He writes (D'Arcy Thompson's translation)\footnote{\emph{History of animals} \cite{A-HA} 487a.}: 
  \begin{quote}\small Of the substances that are composed of parts uniform (or homogeneous)
with themselves, some are soft and moist, others are dry and solid.
The soft and moist are such either absolutely or so long as they are
in their natural conditions, as, for instance, blood, serum, lard,
suet, marrow, sperm, gall, milk in such as have it flesh and the like;
and also, in a different way, the superfluities, as phlegm and the
excretions of the belly and the bladder. The dry and solid are such
as sinew, skin, vein, hair, bone, gristle, nail, horn (a term which
as applied to the part involves an ambiguity, since the whole also
by virtue of its form is designated horn), and such parts as present
an analogy to these. 
\end{quote} 

Thom,\index{Thom, René (1923--2002)} in Chapter 7 of his \emph{Esquisse}, comments on Aristotle's\index{Aristotle (384--322 BCE)} methods of classification developed in the \emph{Parts of animals}, relating it to his own work as a topologist. He recalls that Aristotle 
\begin{quote}\small
attacks therein the Platonic method of \emph{dichotomy}, suggesting in its place\index{place} an interrogative method for taking the substrate into consideration. Thus, if we propose to attain a definition characterizing the ``essence" of an animal, we should not, says Aristotle\index{Aristotle (384--322 BCE)}, pose series of questions bearing on ``functionally independent" characteristics. For example, ``Is is a winged or a terrestrial animal?" or ``Is it a wild animal or a tame one?" Such a battery of questions bearing on semantic fields ---genera--- unrelated one to the other, can be used in an arbitrary order. The questionnaire may indeed lead to a definition of sorts, but it will be a purely artificial one. It would be more rational to have a questionnaire with a \emph{tree} structure, its ramification corresponding to the substrate. For instance, after the question: ``Is the animal terrestrial?", if the answer is yes, we will ask, ``Does the animal have legs?" If the answer is again yes, we then ask, ``Is the foot all in one (solid), or cloven, or does it bear digits?" Thus we will reach a definition which is at the same time a description of the organism in question. Whence a better grasp of its essence in its phenomenal aspect. Aristotle\index{Aristotle (384--322 BCE)} observes, for instance, that if one poses a dilemma bearing on a private opposition, presence of $A$, absence of $A$, the natural posterity of the absence of $A$ in the question-tree is empty. In a way the tree of this questionnaire is the reflection of a dynamic inside the substrate. It is the dynamic of the blowing-up of the centre of the body (the soul), unique in potentiality, into a multitude of part souls in actuality. In a model of the catastrophe type, it is the ``unfolding" dynamics.
\end{quote}

It may be fitting to cite also the following passage from the \emph{History of animals}, on the mysteries of embryology,\index{embryology} and more precisely on the relation between the local and the global\index{biology}\index{biology!local and global in}\index{local and global!in biology} that plays a central role in this domain:\footnote{\emph{History of animals} \cite{A-HA} 589b30 ff.}
\begin{quote}\small
 The fact is that animals,
if they be subjected to a modification in minute organs, are liable
to immense modifications in their general configuration. This phenomenon
may be observed in the case of gelded animals: only a minute organ
of the animal is mutilated, and the creature passes from the male
to the female form. We may infer, then, that if in the primary conformation
of the embryo an infinitesimally minute but absolutely essential organ
sustain a change of magnitude one way or the other, the animal will
in one case turn to male and in the other to female; and also that,
if the said organ be obliterated altogether, the animal will be of
neither one sex nor the other. And so by the occurrence of modification
in minute organs it comes to pass that one animal is terrestrial and
another aquatic, in both senses of these terms. And, again, some animals
are amphibious whilst other animals are not amphibious, owing to the
circumstance that in their conformation while in the embryonic condition
there got intermixed into them some portion of the matter of which
their subsequent food is constituted; for, as was said above, what
is in conformity with nature is to every single animal pleasant and
agreeable. 
\end{quote} 
 
To close this section, I would like to say a few more words on Aristotle\index{Aristotle (384--322 BCE)} as a scientist, from Thom's point of view, and in particular, regarding the (naive) opposition that is usually made between Aristotelian and Galilean science,\index{Galilean physics} claiming that the mathematization of nature, as well as the so-called ``experimental method"  started with Galileo and other moderns, and not before, and in any case, not with Aristotle. 

The lack of ``mathematization" in Aristotle's\index{Aristotle (384--322 BCE)} physics, together with the absence of an ``experimental method" are due in great part due to the lack of measurements, and this was intentional. In fact, Aristotle (like Thom after him) was interested in the qualitative aspects, and not the quantitative. Furthermore, he\index{Aristotle (384--322 BCE)} was reluctant to think in terms of a ``useful science": he was the kind of scholar who was satisfied with doing science for the pleasure of the intellect. Thom's\index{Thom, René (1923--2002)} vision of Aristotelian science was completely different of the commonly accepted one, and the fact that Aristotle was a proponent of the qualitative vs. the quantitative was in line with his own conception of science.  In his article \emph{Aristote et l'avènement de la science moderne~: la rupture galiléenne} (Aristotle\index{Aristotle (384--322 BCE)} and the advent of modern science: the Galilean break),\index{Galilean physics} published in 1991, he writes \cite[p.\,489]{Avenement}: 
          \begin{quote}\small 
          I would be tempted to say that one can see frequently enough a somewhat paternalistic attitude of condescension regarding Aristotle\index{Aristotle (384--322 BCE)} in the mouth of contemporary scholars. I think that this attitude is not justified. It is usually claimed that the Galilean\index{Galilean physics} epistemological break brought to science a radical progress, annihilating the fundamental concepts of Aristotelian physics. But here too, one must rather see the effect of this brutal transformation as a scientific revolution in the sense of Kuhn,\index{Kuhn, Thomas Samuel (1922--1996)}  that is, a change in paradigms, where the problems solved by the ancient theory stop being objects of interest and disappear from the speculative landscape. The new theory produces new problem, which it can solve more or less happily, but above all, it leads to an occultation of all the ancient problematic which nevertheless continues its underground journey under the clothing of new techniques and new formalisms.
\\
I think that very precisely, maybe since ten years, one can see in modern science the reappearance of a certain number of themes and methods that are close to the Aristotelian doctrine, and I personally welcome this kind of resurgence that I will try to describe for you. What I will talk to you about may not belong as much to Aristotle\index{Aristotle (384--322 BCE)}, of whom, by the way, I have a poor knowledge, than to certain recent evolutions of science that recall, I hope without excessive optimism, certain fundamental ideas of Aristotelian physics and metaphysics.\footnote{Je serais tenté de dire qu'une attitude de condescendance un peu paternaliste se remarque assez fréquemment dans la bouche des savants contemporains à l'égard d'Aristote. Et je pense que cette attitude n'est pas justifiée. Il est courant de dire que la rupture épistémologique galiléenne a amené en science un progrès radical, réduisant à néant les concepts fondamentaux de la physique aristotélicienne. Mais là aussi, il faut voir l'effet de cette transformation brutale plutôt à la manière d'une révolution scientifique au sens de Kuhn, c'est-à-dire comme un changement de paradigmes~: les problèmes résolus par l'ancienne théorie cessant d'être objet d'intérêt et disparaissant du champ spéculatif. La nouvelle théorie dégage des problèmes neufs, qu'elle peut résoudre avec plus ou moins de bonheur, mais surtout elle conduit à occulter toute la problématique ancienne qui n'en poursuit pas moins son cheminement souterrain sous l'habillement des nouvelles techniques et des nouveaux formalismes.
  \\
  Je pense que très précisément, depuis peut-être une dizaine d'années, on voit réapparaître dans la science moderne un certain nombre de thèmes et de méthodes proches de la doctrine aristotélicien, et je salue quant à moi cette espèce de résurgence que je vais essayer de vous décrire. Ce dont je vous parlerai ce n'est donc peut-être pas tant d'Aristote, que je connais mal d'ailleurs, que de certaines évolutions récentes de la science qui me semblent évoquer, sans optimisme excessif, j'espère, certaines idées fondamentales de la physique et de la métaphysique aristotélicienne.}
  \end{quote}
In his development, Thom\index{Thom, René (1923--2002)} talks about Aristotelian logic,\index{Aristotelian logic} closely linked to his ontology, and makes science appear as a ``logical language", forming an ``isomorphic image of natural behavior." He also talks about the ``mathematization of nature" that the so-called Galilean break\index{Galilean physics} brought as a hiatus between the mathematical and the common languages, and about the disinterest in the Aristotelian notion of formal causality that characterizes modern science and which would have been so useful in embryology.\index{embryology} He gives examples from Aristotelean mechanics,\index{Aristotelean mechanics} in particular his concept of time, which he links to ideas on thermodynamics, entropy and Boltzmann's\index{Boltzmann, Ludwig (1844-1906)} \index{Boltzmann, Ludwig (1844-1906)!\emph{H}-Theorem} which describes the tendency of an isolated ideal gaz system towards a thermodynamical equilibrium state.

Even if Thom's point of view is debatable, it has the advantage of making history richer and gives it a due complexity. 

As a final remark on mathematics in Aristotle, I would like to recall his reluctance to accept the Pythagorean theories of symbolism in numbers, which he expressed at several places in his work.\footnote{Cf. for instance \emph{Metaphysics} \cite{A-Metaphysics} 1080b16-22.} Incidentally, mathematicians tend to refer  to Plato\index{Plato (428--348 BCE)} rather than to Aristotle\index{Aristotle (384--322 BCE)}, because the former's conception of the world gives a more prominent place to mathematics. But whereas Plato's mathematics is an abstract mathematics, dissociated from the world, the one of Aristotle\index{Aristotle (384--322 BCE)}, is connected with nature. This is consistent with Thom's view. In his paper \emph{Logos phenix} published in the book \emph{Modèles mathématiques de la morphogenèse} \cite{Morpho} which we already mentioned, he writes: ``What remains in me of a professional mathematician can hardly accept that mathematics is only a pointless construction without any attach to reality."\footnote{Ce qui reste en moi du mathématicien professionnel admet difficilement que la mathématique ne soit qu'une construction gratuite dépourvue de toute attache au réel.}

We shall talk more about Thom's\index{Thom, René (1923--2002)} view on Aristotle in the next sections.
  \section{Thom on Aristotle}\label{s:Thom}

Chapter 5 of Thom's \emph{Esquisse} is titled \emph{The general plan of animal organization} and is in the lineage of the zoological treatises of Aristotle\index{Aristotle (384--322 BCE)}, expressed in the  language of modern topology.\index{topology} Thom\index{Thom, René (1923--2002)} writes in the introduction:  
``This presentation might be called an essay in transcendental anatomy, by which I mean that animal organization will be considered here only from the topologists' abstract point of view." He then writes: 
``We shall be concerned with ideal animals, stylized images of existing animals, leaving aside all considerations of quantitative size and biochemical composition, to retain only those inter-organic relations that have a topological and functional character."  In  \S B of the same chapter, Thom returns to Aristotle's\index{Aristotle (384--322 BCE)} notion of homeomerous\index{homeomerous} and anhomeomerous,\index{anhomeomerous} in relation with the stratification\index{stratification!of an animal organism} of an animal's organism, formulating in a modern topological language the condition for two organism to ``have the same organisation."

An organism, in Thom's words, is a three-dimensional ball $O$ equipped with a stratification,\index{stratification!of an animal organism} which is finite if we decide to neglect too fine details. For example, when considering the vascular system, Thom\index{Thom, René (1923--2002)} stops at the details which may be seen with the naked eye: arterioles and veinlets. He writes: ``We will thus avoid introducing fractal morphologies which would take us outside the mathematical schema of stratification."  Seen from this point of view, the homeomerous\index{homeomerous} parts are the strata of dimension three (blood, flesh, the inside of bones...), the two-dimensional strata are the membranes: skin, mucous membrane, periosteum, intestinal wall, walls of the blood vessels, articulation surfaces, etc., the one-dimensional strata are the nerves:  vessel axes, hairs, etc., and the zero-dimensional strata are the points of junction between the one-dimensional strata or the punctual singularities: corners of the lips, ends of hairs, etc. 

Thom\index{Thom, René (1923--2002)} says that two organisms $O$ and $O'$ ``have the same organisation"  if there exists a homeomorphism $h:O\to O'$  preserving this stratification.\index{stratification} He claims that such a formalism generalizes  D'Arcy Thompson's\index{Thompson, D'Arcy Wentworth (1860-1948)} famous diagrams and makes them more precise. 
 The reference here is to Thompson's sketches  from \emph{On growth and Form} which describe passages between various species of fishes using homeomorphisms that preserve zero-dimensional,  one-dimensional, two-dimensional and three-dimensional strata.  
\begin{center}
\includegraphics[width=.8\linewidth]{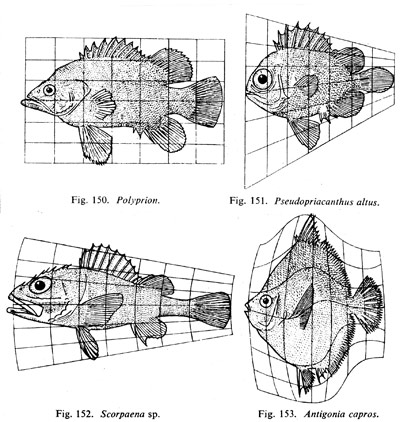} 
\vskip .1in \label{DArcyThom2}
 {\small Homeomorphisms preserving stratifications, from  D'Arcy Thompson's \emph{On growth and form}} 
  \end{center}
  \vskip .3in 
    
Thom\index{Thom, René (1923--2002)} developed his ideas on the stratification of the animal body\index{stratification!of an animal organism}  in the series of lectures given in 1988 at the Solignac Abbey, a medieval monastery in the Limousin (South of France). The lectures are titled \emph{Structure et fonction en biologie aristotélicienne} (\emph{Structure and function is Aristotelian biology}) and the lecture notes ara available  \cite{A-Structure}. 
On p. 7 of these notes, he addresses the question of when two animals have isomorphic stratifications, and he uses for this the notion of isotopy between stratified spaces: two sets $E_1, E_2$ have isotopic stratifications if there exists a stratification\index{stratification} of the product $E\times [0,1]$ such that the canonical projection $p:E\times [0,1]\to[0,1]$ is of rank one on every stratum of $E$, with $E_1=p^{-1}(1)$ and $E_2=p^{-1}(2)$. He considers that this notion is implicitly used by Aristotle in his classification of the animals, insisting on the fact that the latter neglected all the quantitative differences and was only interested in the qualitative ones. In the same passage, he recalls that D'Arcy Thompson, who translated Aristotle's\index{Aristotle (384--322 BCE)} \emph{History of  Animals},  acknowledged that he found there the idea of his diagrams.

        Chapter 7 of Thom's \emph{Esquisse} is called  \emph{Perspectives in Aristotelian biology}. The first part concerns topology\index{topology} and bears the title \emph{The primordial topological intuitions of Aristotelianism: Aristotle and the continuum}. It starts as follows \cite[p.\,165]{Esquisse}: 
        \begin{quote}\small
        We shall present here those intuitions which we believe sub-tend all Aristotelianism. They are ideas that are never explicitly developed by the author, but which---to my mind---are the framework of the whole architecture of his system. We come across these ideas formulated ``by the way" as it were, condensed into a few small sentences that light up the whole corpus with their bright concision.
        \end{quote}
          Thom\index{Thom, René (1923--2002)} highlights and  comments on several citations from Aristotle\index{Aristotle (384--322 BCE)} which show, according to him, that Aristotle\index{Aristotle (384--322 BCE)} was aware of the basic notions of topology.\index{topology} He declares in particular that a careful reading of Aristotle's \emph{Physics}\index{Aristotle (384--322 BCE)!\emph{Physics}} shows that the Philosopher understood the topological distinction between a closed and an open set. He writes, \cite[p.\,167-168]{Esquisse}:           \begin{quote}\smaller 
          Careful reading of the \emph{Physica} leaves little doubt but that [Aristotle] had indeed perceived this difference. ``It is a whole and limited; not, however, by itself, but by something other than itself"\footnote{\emph{Physics} \cite{A-Physics} 207a24-35.} could hardly be interpreted except in terms of a bounded open set. In the same vein, the affirmation: ``The extremities of a body and of its envelope are the same"\footnote{\emph{Physics} \cite{A-Physics} 211b12.} can be identified, if the envelope is of a negligible thickness, with the well-known axiom of general topology: ``Closure of closure is closure itself" expressed by Kuratowski at the beginning of this century. This allows the Stagirite to distinguish two infinites: the great infinite that envelops everything and the small infinite that is bounded. This latter is the infinite of the continuum, able to take an infinity\index{infinity} of divisions (into parts that are themselves continuous).\index{continuous} Whence the definition he proposes: ``The infinite has an intrinsic substrate, the sensible continuum.\index{Aristotle (384--322 BCE)!\emph{Physics}}"\footnote{\emph{Physics} \cite{A-Physics} 208a.} 
            \end{quote}

 In his 1988 article \emph{Les intuitions topologiques primordiales de l'aristot\'elisme} \cite{Thom-Intuitions} and in his 1991 article \emph{Matière, forme et catastrophes} \cite{A-Matiere}, Thom\index{Thom, René (1923--2002)} returns to these matters, explaining that the modern topological distinction between an open and a closed set is related to Aristotle's\index{Aristotle (384--322 BCE)}  philosophical distinctions between form\index{form} and matter\index{form}\index{form!matter and}\index{matter!form and} and between actuality and potentiality.\index{potentiality vs. action} He gives an explanation for the difference between the notion of bounded and unbounded open set in Aristotle's philosophical system: the former may exist as a substrate of being whereas the latter cannot  \cite[p.\,396]{Thom-Intuitions}.
 Concerning the notion of boundary,\index{boundary} he comments on formulae such as: ``Form is the boundary of matter," \cite[p.\,398]{Thom-Intuitions} and ``Actuality is the boundary of potentiality" (\cite[p.\,399]{Thom-Intuitions} and \cite[p.\,380]{A-Matiere}). He recalls that the paradigmatic substance for Aristotle\index{Aristotle (384--322 BCE)} is the living being, which is nothing else than a ball in Euclidean space, whose boundary is a sphere (provided, Thom\index{Thom, René (1923--2002)} says, one neglects the necessary physiological orifices), that is, a closed surface \emph{without boundary}. Shapeless matter is enveloped by form---\emph{eidos}---in the same way as the boundary of a bronze statue defines its shape. The boundary\index{boundary!of an organism} of a living organism is its skin, and its ``interior" exists only as a potentiality. 
 A homeomerous\index{homeomerous} part of an animal has generally a boundary structure constituted by anhomeomerous\index{anhomeomerous} parts. Thus, the substrate of a homeomerous part is not closed.\footnote{Une partie homéomère a en général un bord constitué d'anhoméomères ; ainsi le substrat d'un homéomère n'est pas -- en général -- un ensemble fermé au sens de la topologie moderne. (\cite[p.\,398]{Thom-Intuitions}).} Mathematics, philosophy and biology\index{biology} are intermingled in this interpretation of a living organism, following Aristotle, with formulae like: 
``The opposition homeomerous-anhomeomerous is a ``representation" (a homomorphic image) of the metaphysical opposition: potentiality-act.\index{potentiality vs. action} As the anhomeomerous\index{anhomeomerous} is part of the boundary\index{boundary} of a homeomerous\index{homeomerous} of one dimension higher, we recover a case of the application of \emph{act as boundary of the potentiality}."\footnote{L'opposition homéomère-anhoméomère est une ``représentation" (une image homomorphe) de l'opposition métaphysique : puissance-acte. Comme l'anhoméomère est partie du bord d'un homéomère de dimension plus grande, on retrouve ainsi un cas d'application de l'\emph{acte bord de la puissance} \cite[p.\,400]{Thom-Intuitions}.}
 
  Thom\index{Thom, René (1923--2002)} refers to Chapter 16 of Book Z of the \emph{Metaphysics}\index{Aristotle (384--322 BCE)!\emph{Metaphysics}}, interpreting a sentence of the \emph{Metaphysics}\index{Aristotle (384--322 BCE)!\emph{Metaphysics}}  on entelechy (\tg{>entel'eqeia}) (the ``creative principle", by which being passes from potentiality to action):\index{potentiality vs. action} ``Entelechy separates", as the opposition between an open and a closed set, more precisely, by the fact that ``a closed segment is `actual' as opposed to the semi-open interval which exists only as a\index{potentiality vs. action} potentiality."\footnote{Le segment fermé est ``actuel" par opposition à l'intervalle semi-ouvert qui n'est qu'en puissance.} He writes again, in the same passage  \cite[p.\,381]{A-Matiere}: ``An open substrate characterizes potential entity. A closed substrate is required for the acting being."\footnote{Un substrat ouvert caractérise l'entité en puissance. Un substrat fermé est requis pour l'être en acte.}

In a passage of the \emph{Esquisse}, Thom talks about Aristotle as ``the philosopher of the continuous"\index{continuous} \cite[p.\,viii]{Esquisse}, and he considers that his chief merit was that he was ``the only one  who thought in terms of the continuous" for hundreds, may be thousands of years.
 In his article \emph{Logos phenix}, Thom\index{Thom, René (1923--2002)} writes: ``How can we explain that mathematics represents the real? The response, I thinks, is offered to us by the intuition of the continuous.\index{continuous}  [\ldots] The introduction of an underlying continuous substrate allows us henceforth to explain the significant---non-trivial---character of several mathematical theorems."\footnote{Comment expliquer que les mathématiques représentent le réel ? La réponse, je crois, nous est offerte par l'intuition du continu. [\ldots] L'introduction d'un substrat continu sous-jacent permet dès lors de s'expliquer le caractère signifiant -- non trivial -- de bien des théorèmes de la mathématiques.} \cite[p.\,292ff]{Morpho}. Thom expands then this idea, giving as examples the implicit function theorem and Taylor's formula, each of them allowing a global knowledge from a local one, all this based on the existence of a continuous substrate.\index{continuous} In the book \emph{Prédire n'est pas expliquer} \cite[p.\,81-82]{A-Predire}, turning to the distinction that the Greek philosophers made between the discrete and the continuous,\index{continuous!vs. discrete}\index{discrete!vs. continuous} he declares: ``For me, the fundamental aporia of mathematics is certainly the opposition between the discrete and the continuous.\index{continuous!vs. discrete}\index{discrete!vs. continuous} And this aporia at the same time dominates all thought. [\ldots] The origin of scientific thinking, we find it in the apories of Zeno of Elea:\index{Zeno of Elea (c. 490--430 BCE)}\index{Zeno of Elea (c. 490--430 BCE)!paradoxes} the story of Achilles and the tortoise. Here, we find the crucial opposition between the discontinuous and the continuous."\footnote{Pour moi, l'aporie fondamentale de la math\'ematique est bien dans l'opposition discret-continu. Et cette aporie domine en même temps toute la pens\'ee. [\ldots] L'origine de la pens\'ee scientifique, on la trouve dans les apories de Z\'enon l'\'El\'ee : l'histoire d'Achille et de la tortue. Il y a l\`a l'opposition cruciale entre discontinu et continu.}   Thom\index{Thom, René (1923--2002)} declares in the same work, concerning this distinction between the continuous\index{continuous} and the discrete:\index{continuous!vs. discrete}\index{discrete!vs. continuous}\index{discrete}  
  \begin{quote}\small
  The continuous\index{continuous} is in some way the universal substrate of thought, and in particular of mathematical thought. But we cannot think of anything in an effective way without having something like the discrete in the continuous\index{continuous!vs. discrete}\index{discrete!vs. continuous} flow of mental processes: there are words, sentences, etc. The \emph{logos},  discourse, is always discrete; these are words, coming in with a certain order, but they are discrete words. And the discrete\index{discrete} immediately calls down the quantitative. There are points: we can count them; there are words in a sentence: we can classify them quantitatively by the grammatical function they occupy in a sentence. However there is an undeniable multiplicity."\footnote{Le continu est en quelque sorte le substrat universel de la pens\'ee, et de la pens\'ee math\'ematique en particulier. Mais on en peut rien penser de mani\`ere effective sans avoir quelque chose comme le discret dans ce d\'eroulement continu de processus mentaux : il y a des mots, il y a des phrases, etc. Le \emph{logos}, le discours, c'est toujours du discret ; ce sont des mots entrant dans une certaine succession, mais des mots discrets. Et le discret appelle imm\'ediatement le quantitatif. Il y a des points : on les compte ; il a des mots dans une phrase : on peut les classer quantitativement par la fonction grammaticale qu'ils occupent dans la phrase, mais il n'empêche qu'il y a une incontestable multiplicit\'e.}
  \end{quote} 
  On the same subject, in his article \emph{Logos phenix}, Thom\index{Thom, René (1923--2002)} writes \cite[p.\,294]{Morpho}:
  \begin{quote}\small Meaning is always tied to the attribution of a place of spatial nature to an expression formally encoded. There should always be, in any meaningful message, a discontinuous component tied to the generative mechanisms of language---to symbols---, and a continuous component, a substrate, in which the continuous component cuts out a place.\footnote{Le sens est toujours lié à l'attribution d'une place de nature spatiale à une expression formellement codée. Il y aurait toujours, dans tout message signifiant, une composante discontinue liée aux mécanismes génératifs de la langue -- aux symboles --, et une composante continue, un substrat, dans lequel la composante continue découpe une place.}
  \end{quote}

In the foreword to the \emph{Esquisse}, Thom\index{Thom, René (1923--2002)} writes that Aristotle's\index{Aristotle (384--322 BCE)} geometrical insight was founded uniquely on an intuition of the continuous, where a segment of a straight line is not made out of points but of sub-segments.  Neither Dedekind nor Cantor, he says,\index{Thom, René (1923--2002)} have taken that road. He writes that the single isolated point (the one we consider when we take a point $O$ on the $x'x$-axis), exists only ``potentially." Such a point, according to him,  aspires to actuality by duplicating itself into two points $O_1, O_2$, $O_1$ adhering to the left, $O_2$ to the right. ``These two points then being \emph{distinct} even though they are \emph{together} (\tg{<'ama}), the two half-segments so limited attain full existence, being in actuality."  \cite[p.\,viii]{Esquisse} Thom\index{Thom, René (1923--2002)} refers for that to the \emph{Metaphysics}\index{Aristotle (384--322 BCE)!\emph{Metaphysics}} 139a3-7\index{Aristotle (384--322 BCE)!\emph{Physics}}.\footnote{In a footnote, Thom refers to a passage in Dieudonné's\index{Dieudonné, Jean (1906--1992)} \emph{Pour l'Honneur de l'esprit humain} \cite{Dieudonne} as a ``fine example of modern incomprehension of the Aristotelian point of view." In p. 229 of this book, Dieudonné criticizes Aristotle\index{Aristotle (384--322 BCE)} for his view of infinity,\index{infinity} on the basis of the passage 231a21-232a22 of the \emph{Physics} where he discusses points on a straight line, describing his reasoning as an example of ``mental coonfusion."}
The question is also that of knowing what expressions like ``is made of" or ``are part of", applied to points and lines, mean. In his Solignac 1988 notes, quoting a passage from Aristotle's\index{Aristotle (384--322 BCE)} \emph{Parts of animals} in which the latter compares an animal vascular system  to a garden irrigation ditch,\footnote{\emph{Parts of animals} \cite{A-Parts} 668a10-13.} Thom\index{Thom, René (1923--2002)} notes that the Philosopher  goes as far as to say that blood is not part of the organism,\footnote{\emph{Parts of animals} \cite{A-Parts} 636b21.} because it is dense there, and that this does not belong to the definition of ``being part of" \cite[p.\,7]{A-Structure}.

   In the same foreword to his \emph{Esquisse}, Thom\index{Thom, René (1923--2002)} recalls that Aristotle\index{Aristotle (384--322 BCE)} considered that the underlying substrate of both matter and form\index{form}\index{form!matter and}\index{matter!form and} is continuous\index{continuous} \cite[p.\,{\sc viii}]{Esquisse}:
    \begin{quote}\small
    I knew of course that the hylomorphic schema---of which I make use in the catastrophe formalism---originated in the Stagirite's work. But I was unaware of the essential fact that Aristotle\index{Aristotle (384--322 BCE)} had attempted in his \emph{Physics} to construct a world theory based not on numbers but on continuity\index{continuous}. He had thus (at least partly) realized something I have always dreamed of doing---the development of mathematics of the continuous,\index{continuous} which would take the notion of continuum as point of departure, without (if possible) any evocation of the intrinsic generativity of numbers.
    \end{quote}
  
On the same subject, Thom\index{Thom, René (1923--2002)} recalls in the \emph{Esquisse} that Aristotle's decision to quit Plato's\index{Plato (428--348 BCE)} Academy is due to a disagreement with his master concerning the notion of continuity\index{continuous}.  He explains this in a long passage \cite[p.\,166]{Esquisse}:  
\begin{quote}\small
The Ancients knew that the moving point generates a curve, that a moving curve generates a surfaces, and that the movement of a surface generates a volume. It seems that the aging Plato---or his epigones---considered this generation to be of the type of discrete generativity, that of the sequence of natural integers. So the point, which is a pure ``zero", could not serve as a base of this construction---whence the necessity of ``thickening" the point into a ``unsecable length"\index{unsecable length}\index{indivisible line}  (\tg{>'atomos gramm'h}), which was the generating principle of the straight line (\tg{>arq`h gramm\~hs}). The \emph{Timaeus}'\index{Plato (428--348 BCE)!\emph{Timaeus}} demiurge could then use this unsecable length\index{indivisible line} to construct the polygons and polyhedrons which constitute the elements.  It is odd to note that this kind of hypothesis still haunts our contemporary physicists; the elementary length ($10^{-33}$ cm) below which space no longer has a physical meaning, or that absolute spacial dimension given by the confinement of quarks in nuclear physics, are so many absolute ``lengths" associated with physical agents. Why did Aristotle\index{Aristotle (384--322 BCE)} reject this sort of hypothesis? No doubt because he held number generativity in disregard. His revolt against Plato\index{Plato (428--348 BCE)} is that of the topologist against the arithmetic imperialism, that of the apostle of the qualitative against the quantitative. Aristotle\index{Aristotle (384--322 BCE)} basically postulates the notion of \emph{continuity\index{continuous}} (\tg{suneq'es}), and it is in the name of the divisibility of the continuum that he refuses the ``indivisible lines."\index{indivisible line}\index{unsecable length}\index{indivisible line} \emph{A priori} this is a paradoxical position. Indeed, Aristotle never admitted the existence of space in the sense in which we have considered it since Descartes.\index{Descartes, René (1596--1650)} We know why: his substantialist metaphysics required that extent is made a predicate of the substance (the \emph{topos}); in no way could substance, matter, be a predicate of space. For Aristotle space is never generated by some intrinsic generative mechanism as our Cartesian space is generated by the $\mathbb{R}^3$ additive group if translations; at the most it is the place of some entity (\emph{ousia}), for it is never empty. This decision to relegate space to a kind of total ostracism led him, by a singular rebound, to multiply the kinds of matter. Each time of change (\tg{metabol'h}), each genus (\tg{g'enos}), requires a specific matter. But all these matters have one thing in common: they are continua (\tg{suneq'es}); in this sense they all have the character of spatial extension.
\end{quote}
 
 Needless to say, the question of the discreteness or continuity\index{continuous!vs. discrete}\index{discrete!vs. continuous} of the ultimate constitution of nature has been at the edge of philosophical thought, at least since the 5$^{\mathrm{th}}$ century BC; one thinks of Leucippus, Democritus, Empedocles and their followers. In the modern period, as mathematicians, we think of Riemann\index{Riemann, Bernhard (1826--1866)} who addressed the question in his 
 Habilitation lecture \emph{\"Uber die Hypothesen, welche der Geometrie zu Grunde
liegen} (On the hypotheses that lie at the bases of geometry) \cite{Riemann-Ueber} (1854) (see also the discussion in \emph{Riemann1}), and Hermann Weyl,\index{Weyl, Hermann (1885--1955)} who continued Riemann's tradition. 
  
 Plato\index{Plato (428--348 BCE)} disliked the notion of point, considering that it has no geometric meaning. To compensate for this, he used the notion of ``unsecable", that is, ``indivisible" length.\index{unsecable length}\index{indivisible line} Aristotle\index{Aristotle (384--322 BCE)}, in a passage of the \emph{Topics}\index{Aristotle (384--322 BCE)!\emph{Topics}},\footnote{\emph{Topics}, \cite{A-Topics} 148b27.} quotes Plato's definition of straight lines. At several places in his writings, he discussed the relation between a point and a line.  For him, a line is a continuous object and as such it cannot be neither a collection of points nor a collection of indivisible objects of any kind.  A point may only be the start, or the end, or a division point of a line, but it is not a magnitude. The question of the nonexistence of indivisibles is so important for him  that it is treated in several of his writings.\index{Aristotle (384--322 BCE)!\emph{On the heavens}}\index{Aristotle (384--322 BCE)!\emph{Physics}}\footnote{\emph{Physics} \cite{A-Physics} 215b12-22,  220a1-21 et  231b6, \emph{On generation and corruption}\index{Aristotle (384--322 BCE)!\emph{On generation and corruption}}  \cite{A-Generation} 317a11, \emph{On the heavens} \cite{A-Heavens} III. 1, 299a10 ff.; there are other passages.}  In his treatise \emph{On the heavens},\index{Aristotle (384--322 BCE)!\emph{On the heavens}} Book~III, Chapter 1, he discusses the analogous question in higher dimensions, that is, the impossibility of a surface to be a collection of lines, and of a solid to be a collection of surfaces, unless, he says, we change the axioms of mathematics, and he adds that this is not advisable. Talking of a change in the axioms in order to obtain a result is a typical attitude of Aristotle acting as a mathematician. Entering into the question of whether a geometric line, or, more generally, a geometric body, is constituted of its points, and if yes, in what sense this is so,   leads us deep into considerations which several Greek philosophers have thoroughly considered (we mentioned Plato\index{Plato (428--348 BCE)} and Aristotle\index{Aristotle (384--322 BCE)}, but these questions were extensively studied before them, especially by Zeno of Elea).\index{Zeno of Elea (c. 490--430 BCE)} One may think of the axiomatization of the real line realized in the 19$^{\mathrm{th}}$ century, by Cantor, Dedekind and others.  Let me quote here a sentence by Plotnitsky in the present volume \cite{Plotnitsky} which expresses the current veiwpoint: ``In sum, we do not, and even cannot, know how a continuous line, straight or curved (which does not matter
topologically), is spatially constituted by its points, but we have algebra to address this question,
and have a proof that the answer is rigorously undecidable."

 It might be recalled here that there exists a treatise called \emph{On indivisible lines}\index{indivisible line}\index{unsecable length}\index{indivisible line} and belonging to the Peripatetic school (may be to Aristotle), in which the author criticizes item by item the arguments of the disciples of indivisible lines\index{indivisible line}\index{Aristotle (384--322 BCE)!\emph{On indivisible lines}} \cite{A-Indivisible}.

On the relation between Plato\index{Plato (428--348 BCE)} and Aristotle\index{Aristotle (384--322 BCE)}, Thom\index{Thom, René (1923--2002)} writes the following, in a note on p. 186 of his \emph{Esquisse}: ``The relations between Plato and Aristotle\index{Aristotle (384--322 BCE)} constitute one of the \emph{topoi} of philosophical erudites. [\ldots] My own position on the question is that of an autodidact."\footnote{Thom  refers to the books \cite{Robin} by Robin and \cite{Cherniss} by Cherniss. Concerning this subject, I would like to refer the reader to the appendix to my article, written by S. Negrepontis.} For more on this subject, the reader is referred to the appendix to the present article.

Let us return to the notion of boundary.\index{boundary}

It is interesting to see that this notion is included in Euclid's \emph{Elements}\index{Euclid (c. 325--c. 270 BCE)!\emph{Elements}} among the elementary notions, at the same level as ``point", ``line", ``angle," etc.  The boundary, there, is what defines a figure. Definition 14 of Book I says: ``A figure is that which is contained by any boundary or boundaries."\index{boundary} This is clearly related to Thom's\index{Thom, René (1923--2002)} idea, following Aristotle, that a form\index{form} is defined by its boundary.   Euclid\index{Euclid (c. 325--c. 270 BCE)} also uses the notion of boundary\index{boundary} when he talks about the measure of angles (not only rectilinear angles).  Proclus,\index{Proclus of Lycia (411--485)} in his \emph{Commentary on the First Book of Euclid's Elements},\index{Proclus of Lycia (411--485)!\emph{In Euclidem}}  writes that the notion of boundary belongs to the origin of geometry since this science originated in the need to measure areas of pieces of land. Aristotle\index{Aristotle (384--322 BCE)}, in the \emph{Physics},\index{Aristotle (384--322 BCE)!\emph{Physics}} says that a body may be defined as being ``bounded by a surface."\footnote{\emph{Physics} \cite{A-Physics} 204b.}

 The notion of form,\index{form} since the origin of geometry, is closely related to the notion of boundary,\index{boundary} and it is not surprising to see that the mathematical notion of boundary which, needless to say, was essential in Thom's\index{Thom, René (1923--2002)} mathematical work, is also central in his philosophical thought.
He declares in the interview \emph{La th\'eorie des catastrophes} conducted in 1992 \cite{INA}: 
\begin{quote}\small
All the unity of my work  is centered at the notion of boundary,\index{boundary} since the notion of cobordism\index{cobordism} is only one of its generalizations. The notion of boundary seems to me the more important today since I am interested in Aristotelian metaphysics. For Aristotle\index{Aristotle (384--322 BCE)}, the boundary\index{boundary} is an individualisation principle. The marble statue is matter, in the block from where the sculptor extracted it, but its is its boundary which defines its form.\index{form}\footnote{Toute l'unit\'e de mon travail tourne autour de la notion de bord, car la notion de cobordisme n'en est qu'une g\'en\'eralisation. La notion de bord me paraît d'autant plus importante aujourd'hui que je m'int\'eresse à la m\'etaphysique aristot\'elicienne. Pour Aristote, le bord est principe d'individuation. La statue de marbre est mati\`ere dans le bloc d'où le sculpteur l'a tir\'ee, mais c'est son bord qui d\'efinit sa forme.}
\end{quote}

  Thom also talks about the unity of his work in the series of interviews \emph{Pr\'edire n'est pas expliquer}. He declares (1991, \cite{A-Predire}, p. 20--21):
    \begin{quote}\small
   Truly, there exists a real unity in my reflections. I can see it only today, after I pondered a lot about it, at the philosophical level. And this unity, I find it in the notion of boundary.\index{boundary} That of cobordism\index{cobordism} is related to it. The notion of boundary is all the more important since I was immersed into Aristotelian metaphysics. For Aristotle\index{Aristotle (384--322 BCE)}, a being, in general, is what is here, separated. It possesses a boundary,\index{boundary} it is separated from the ambiant space. In other words, the boundary of an object is its form.\index{form} A concept has also a boundary,\index{boundary} viz. the definition of that object. On the other hand, this idea that the boundary defines the object is not completely exact for a topologist. It is only true in the usual space. But the fact remains that, starting from this notion of boundary,\index{boundary} I developed a few mathematical theories that were useful to me; then I looked into the applications, that is, on the possibilities of sending a space into another one, in a continuous manner. From here, I was led to study cusps and folds, objects that have a mathematical meaning.\footnote{``En v\'erit\'e, il existe une r\'eelle unit\'e dans ma r\'eflexion. Je ne la per\c cois qu'aujourd'hui, apr\`es y avoir beaucoup r\'efl\'echi, sur le plan philosophique. Et cette unit\'e, je la trouve dans cette notion de bord. Celle de cobordisme lui est li\'ee. [\ldots]
   La notion de bord est d'autant plus importante que j'ai plong\'e dans la m\'etaphysique aristot\'elicienne. Pour Aristote, un être, en g\'en\'eral, c'est ce qui est l\`a, s\'epar\'e. Il poss\`ede un bord, il est s\'epar\'e de l'espace ambiant. En somme, le bord de la chose, c'est sa forme. Le concept, lui aussi, a un bord : c'est la d\'efinition de ce concept. Cette id\'ee que le bord d\'efinit la chose n'est d'ailleurs pas tout \`a fait exacte pour un topologue. Ce n'est vrai que dans l'espace usuel.
  Il reste que, partant de cette notion de bord, j'ai d\'evelopp\'e, quelques th\'eories math\'ematiques qui m'ont servi ; puis je me suis pench\'e sur les applications, c'est-\`a-dire sur les possibilit\'es d'envoyer un espace dans un autre, de mani\`ere continue. J'en ai \'et\'e amen\'e \`a \'etudier les fronces et les plis, objets qui ont une signification math\'ematique." 
     }
     
     \end{quote}

Let us return to the Stokes formula\index{Stokes formula}. It  establishes a precise relation between a domain and its boundary.\index{boundary} In his article  \emph{Aristote topologue}  \cite{Aristote-Topologue} (1999), Thom\index{Thom, René (1923--2002)} writes that one can interpret a passage of Aristotle's \emph{Physics}\index{Aristotle (384--322 BCE)!\emph{Physics}} in which he talks about the ``minimal limit of the enveloping body"\footnote{\emph{Physics} \cite{A-Physics} 211b11.}  as the homological Stokes formula: $d\circ d=0$, the dual of the usual Stokes formula\index{Stokes formula}, concerning differential forms and expressing the fact that the boundary of the boundary\index{boundary} is empty. ``This formula, he writes, essentially expresses the closed character of a human being. Because if there is boundary, then there is blood loss, with a threat to life. Hence the role of the operator $d^2=0$ from homological algebra as an \emph{ontology detector}, and its profound biological\index{biology} interpretation."\footnote{La formule explique essentiellement le caract\`ere clos de l'\^etre vivant. Car s'il y a un bord, il y a perte de sang, avec menace pour la vie. D'où le rôle de \emph{d\'etecteur d'ontologie} qu'est l'op\'erateur bord $d^2=0$ de l'alg\`ebre homologique et de sa profonde interpr\'etation biologique."} This identification of the Stokes formula\index{Stokes formula} with a formula in Aristotle\index{Aristotle (384--322 BCE)} was already carried on in his 1991 article \emph{Matière, forme et catastrophes} \cite[p.\,381]{A-Matiere}.
   In the series of interviews \emph{ Prédire n'est pas expliquer} (1991), he\index{Thom, René (1923--2002)} also talks about the homological Stokes  formula\index{Stokes formula} $d^2=0$ and his biological interpretation: ``The boundary\index{boundary} of the boundary is empty; this is the great axiom of topology\index{topology} and of differential geometry in mathematics, but it is an expression of the spacial integrity of the boundary of the organism."\footnote{Le bord du bord est vide ; c'est le grand axiome de la topologie, de la g\'eom\'etrie diff\'erentielle en math\'ematiques, mais cela exprime l'int\'egrit\'e spatiale du bord de l'organisme.} \cite[p.\,111]{A-Predire} We recall incidentally that the Stokes formula\index{Stokes formula} is the basis tool in the proof of the theorem stating that the characteristic numbers of two cobordant manifolds, computed from the tangent bundles, coincide, and that the great theorem of Thom,\index{Thom, René (1923--2002)} the one for which he was awarded the Fields medal, is the converse of this one.

 We should also talk about the relation between the local and the global in biology,\index{biology}\index{biology!local and global in}\index{local and global!in biology} and in particular, the problems of ``local implies global" type. 
 
 Thom\index{Thom, René (1923--2002)} declares in his 1992 interview on catastrophe theory\index{catastrophe theory} \cite{INA} that the big problem of biology is the relation between the local and the global, that this is a philosophical problem which has to do with ``extent" and that, at the same time, it is the object of topology.\index{topology} Catastrophe theory\index{catastrophe theory} has something to say about the common features of the evolution of the form of a wave, of a cloud, of a living cell, of a fish and of any other living being, but also on the question of how does morphogenesis---the birth of form---affect the later development of a form.\index{form} This is the biological counterpart\index{biology} of the question of how the local implies the global. Thom,\index{Thom, René (1923--2002)} in the interview, recalls that topology is essentially the study of the ways that make a relation between a given local property and a global property to be found, or conversely: given a global property of a space,\index{local and global} to find its local properties, around each point. He concludes by saying that there is a profound methodological unity between topology and biology.\index{biology}  
 
 Talking about Aristotle\index{Aristotle (384--322 BCE)} in a mathematical context, one is tempted to say a few words about the logic he founded, the so-called Aristotelian logic,\index{Aristotelian logic} which is different from the abstract mathematical logic---a 19$^{\mathrm{th}}$ century invention.  Thom had his own ideas on the matter, and, needless to say, if the question of which among the two logics is more suitable to science may be raised, his preference goes to Aristotelian logic. We leave to him the final word of this section, from his article \emph{Aristote et l'avènement de la science moderne}, \cite[p.\,489]{Avenement}:
 
       \begin{quote}\small
       Aristotle\index{Aristotle (384--322 BCE)} had perfectly understood that there is no pertinent logic without an ontology which serves for it as foundations. In other words, if a logic may serve to describe in an efficient way certain aspects of the real world, if logical deduction is a reflexion of the behavior of real phenomena, then, that logic must necessarily be connected with the reality of the external world. And indeed, it is very clear that Aristotle's logic was one with his physics. [\ldots] I will say in most formal way that the so-called progress of logic, realized since the appearance in the 19$^{\mathrm{th}}$  century of formal logic with Boole were in fact Pyrrhic progresses, in this sense that which we gained from the point of view of rigor, we lost it from the point of view of pertinence. Logic wanted to be separate of any ontology, and, for that, it became a gratuitous construction, in some way modeled on mathematics, but such an orientation is even less justified than in the case of mathematics.\footnote{Aristote avait parfaitement compris qu'il n'y a pas de logique pertinente sans une ontologie qui lui sert de fondement. Autrement dit, si une logique peut servir à décrire efficacement  certains aspects du réel, si la déduction logique est un reflet du comportement des phénomènes réels, eh bien, c'est que la logique doit nécessairement avoir un rapport avec la réalité du monde extérieur. Et il est bien clair en effet que la logique d'Aristote faisait corps avec sa physique. [\ldots] Je serai tout à fait formel en disant que les prétendus progrès de la logique, réalisés depuis l'apparition de la logique formelle avec Boole au XIX$^{\mathrm{e}}$ siècle, ont été en fait des progrès à la Pyrrhus, en ce sens que ce qu'on a gagné du point de vue de la rigueur, on l'a perdu du point de vue de la pertinence. La logique a voulu se séparer de toute ontologie et, de ce fait, elle est devenue une construction gratuite, un peu sur le même modèle que les mathématiques, mais une telle orientation est encore moins motivée que dans le cas des mathématiques.}
        \end{quote}

    \section{On form}\index{form}
     Thom\index{Thom, René (1923--2002)} was thoroughly involved in questions of morphogenesis.\index{form} 
   In his article  \emph{Matière, forme et catastrophes}, he recalls that it was in 1978 that for the first time he made the connection between catastrophe theory\index{catastrophe theory} and Aristotle's hylemorphism\index{hylemorphism} theory. In the foreword to his \emph{Esquisse}, he writes \cite[p.\,{\sc viii}]{Esquisse}:  ``I knew of course that the hylomorphic schema---of which I make use in the catastrophe formalism---originated in the Stagirite's work." Aristotle's theory of hylemorphism\index{hylemorphism} is discussed thoroughly in the \emph{Esquisse}. According to that theory, every being (whether it is an object or a living being) is composed in an inseparable way of a matter (\emph{hylé}, \tg{<'ulh})\footnote{The  word  \tg{<'ulh}, before Aristotle\index{Aristotle (384--322 BCE)}, designated  shapeless wood, and the introduction of this word in philosophy is due to Aristotle himself.}  and form (\emph{morphê}, \tg{morf'h}).\index{form} Matter, from this point of view, is a potentiality, a substrate awaiting to receive form\index{form} in order to become a substance---the substance of being, or being itself. 
  In the \emph{Metaphysics}\index{Aristotle (384--322 BCE)!\emph{Metaphysics}}, we can read: ``I call form the essential being and the primary substance of a thing."\footnote{\emph{Metaphysics} \cite{A-Metaphysics} 1032b1-2.} In the treatise \emph{On the soul},\index{Aristotle (384--322 BCE)!\emph{On the soul}} Aristotle\index{Aristotle (384--322 BCE)} states that the soul is the \emph{form} of a human being.\footnote{\emph{On the soul}, \cite{A-Soul} 412a11.} In the \emph{Physics},\index{Aristotle (384--322 BCE)!\emph{Physics}} he writes  that the \tg{f'usis} (nature) of a body is its form,\footnote{\emph{Physics} \cite{A-Physics} 193a30.}
  and that flesh and bone do not exist by nature until they have acquired\index{form} their form.\footnote{\emph{Physics} \cite{A-Physics} 193b.} Matter and form also make the difference between the domain of interest of a physicist and that of a mathematician. The former, according to Aristotle\index{Aristotle (384--322 BCE)}, studies matter and form\index{form}\index{form!matter and}\index{matter!form and},\footnote{\emph{Metaphysics} \cite{A-Metaphysics} 1037a16-17 and \emph{Physics} \cite{A-Physics} 194a15.} whereas the latter is only concerned with form.\index{Aristotle (384--322 BCE)!\emph{Posterior analytics}}\footnote{\emph{Posterior analytics} \cite{A-Posterior},  79a13.} In the \emph{Metaphysics}\index{Aristotle (384--322 BCE)!\emph{Metaphysics}}, he writes that mathematical objects constitute a class of things intermediate between forms and sensibles.\footnote{\emph{Metaphysics} \cite{A-Metaphysics} 1059b9.} In another passage of the \emph{Physics},\index{Aristotle (384--322 BCE)!\emph{Physics}}\footnote{\emph{Physics} \cite{A-Physics}  192a24.} he writes that ``matter desires form as much as a female desires a male." Form, according to him, is \emph{what contains}, and it may even contain the infinite:  ``For the matter and the infinite are contained inside what contains them, while it is form which contains."\footnote{\emph{Physics} \cite{A-Physics}  207a35.} 
  It also appears clearly in these writings that Aristotle\index{Aristotle (384--322 BCE)} did not conceive form as a self-contained entity which lives without matter.  
  
  Thom\index{Thom, René (1923--2002)} completely adhered to Aristotle's theory of form,\index{form} which the latter expanded especially in his biological treatises.\index{biology} Thom highlighted the importance of these ideas in biology,\index{biology} and more particularly in embryology,\index{embryology} namely, the idea of a form\index{form} which tends to its own realization.
             
The paper \emph{L'explication des formes spatiales : réductionnisme ou platonisme} (The meaning of spacial forms: reductionism or platonism) (1980) \cite{Explication} by  Thom concerns the notion of form\index{form} and its classification. Thom\index{Thom, René (1923--2002)} tried there to give a mathematical basis to phenomenal concepts. In this setting, the substrate of a morphology\index{morphology} is the four-dimensional Euclidean space. From a mathematical point of view, a \emph{form} is then a closed subset of space-time\footnote{Thom speaks of space-time in the classical sense, that is, a four-dimensional space whose first three coordinates represent space and the fourth one time. This is not the space-time that is used in the theory of relativity.} up to a certain equivalence relation, and, for him,  one of the fundamental problems in morphology\index{morphology} is to make precise, from the mathematical point of view, this equivalence relation.   In biology,\index{biology} Thom\index{Thom, René (1923--2002)} talks about a ``congruence in the sense of D'Arcy Thompson" \cite{Thompson}. The relation satisfies certain metrical constraints which, he says, ``are generally impossible to formalize." This is, he says, the problem that biometrics has to solve: for instance, to characterize a certain bone of a given animal species. The problem is obviously open, but Thom\index{Thom, René (1923--2002)} adds that often, some subtle psychological mechanisms of form\index{form} recognition will allow one to decide almost immediately whether two objects have the same form\index{form} or not. 
  
   In the introduction to the course his 1980 Solignac course \cite{Solignac}, Thom says that
 biology\index{biology} is a morphological\index{morphology} discipline, concerned with form,\index{form} and topology,\index{topology}  as a branch of mathematics which involves the study of form,\index{form}  is at the basis of theoretical biology.\index{biology} From his point of view, there are two steps in a morphological discipline: the classification---giving names to the various forms, the identification of stable forms, etc.---and, after that, the theorization, namely, building a theory which is ``generative" in the sense that it confers to certain forms (or agregates of form)\index{form} a certain power of determining other forms which are close to them. This program, says Thom,\index{Thom, René (1923--2002)} is partially realized in linguistics,\index{linguistics} which he also considers as a morphological discipline.   In 1971, he published an essay on the subject, \emph{Topologie et linguistique}\footnote{This article was published five years later in Russia, with an introduction by Yuri Manin.} \cite{TL}  in which he develops a general theory of linguistics based on topology\index{topology} and where the accent is on morphology,\index{morphology} again in a pure Aristotelian tradition. Thus, like biology, linguistics\index{linguistics} is, for Thom,\index{Thom, René (1923--2002)} part of the general theory of forms.\index{form} A sentence, a phrase, whether it is written or oral, is, according to him, a form.\index{form} More than that, morphology\index{morphology} is what unifies language---a complicated process whose study pertains at the same time to physiology, psychology, sociology, and other fields.
    
     In his article \emph{Aristote et l'avènement de la science moderne} (1991), discussing the relation between Aristotle\index{Aristotle (384--322 BCE)} and modern science, Thom\index{Thom, René (1923--2002)} writes  \cite[p.\,491ff]{Avenement}: 
             \begin{quote}\small I belong to those who think the the hylemorphic schema is still valid, because it is equivalent to the classifying role of concept in the verbal description of the world. [\ldots] I am convinced that during the last years, in several disciplines, there appeared situations that can be explained by the presence of local fields or forms and that absolutely justify the old Aristotelian hylemorphic model, according to which nature is in some sense captured by form. Of course, I do not hide myself the fact that here, Aristotelian form, the ``\emph{eidos}", was a being that had nothing mathematical. It was an entity that carried its own ``\emph{energeia}", its activity, and it is clear that, tor Aristotle\index{Aristotle (384--322 BCE)}, form\index{form} did not have the status of a mathematical object that would have led him to a certain form of Platonism. The fact remains that the Aristotelian ``\emph{eidos}"  has a certain efficient virtue which, anyway, one has to explain, and in the theories of modern science which I am alluding to, the efficiency of the ``\emph{eidos}" is expressed in mathematical terms, for instance using structural stability.\index{stability}\footnote{Je suis de ceux qui croient que le schème hylémorphique garde toute sa valeur, car il est l'équivalent du rôle classificateur du concept dans la description verbale du monde. [\ldots] Je suis convaincu que ces dernières années ont vu dans un assez grand nombre de disciplines réapparaître des situations qu'on peut expliquer par la présence de champs locaux ou de formes et qui justifient tout à fait à mon avis le vieux modèle hylémorphique d'Aristote, selon lequel la matière en quelque sorte est capturée par la forme. Bien entendu, je ne me dissimule pas qu'ici la forme aristotélicienne, l'``\emph{eidos}", était un être qui n'avait rien de mathématique. C'était une entité qui portait en elle son ``\emph{energeia}", son activité, et il est clair que, pour Aristote, la forme n'avait pas un statut de caractère mathématique qui l'aurait obligé à une certaine forme de platonisme. Mais il n'en demeure pas moins que  l'``\emph{eidos}" aristotélicien a une certaine vertu efficace qu'il faut expliquer de toute façon, et dans les théories de la science moderne auxquelles je fais allusion, l'efficace de  l'``\emph{eidos}" s'exprime en termes mathématiques, par la théorie de la stabilité structurelle par exemple.}
\end{quote} 
 
The book \emph{Structural stability and morphogenesis} starts with a \emph{Program} in which Thom\index{Thom, René (1923--2002)} presents the problem of succession of forms\index{form} as one of the central problems of human thought. He writes: 
\begin{quote}\small
Whatever is the ultimate nature of reality (assuming that this expression has a meaning), it is indisputable that our universe is not chaos.\index{chaos} We perceive beings, objects, things to which we give names. These beings or things are forms or structures endowed with a degree of stability;\index{stability}  they take up some part of space and last for some period of time.  [\ldots] we must concede that the universe we see is a ceaseless creation, evolution, and destruction of forms and that the purpose of science is to foresee this change of form, and, if possible, to explain it.
\end{quote}

  In his article \emph{Aristote et l'avènement de la science moderne}, \cite{Avenement} Thom\index{Thom, René (1923--2002)} declares that since the advent of Galilean physics,\index{Galilean physics} which emphasizes motion in a world in which there is no place for generation and corruption,\index{Aristotle (384--322 BCE)!\emph{On generation and corruption}}\footnote{This is also a reference to Aristotle's \emph{On generation and corruption} \cite{A-Generation}.} considerations on form\index{form} disappeared from physics, even though morphology\index{morphology} is present in biology.\index{biology} Modern science, he says, is characterized by the disappearance of this central notion of form, which played a central role in Aristotle's ontology (\cite[p.\,491]{Avenement}).\footnote{``Le monde de Galilée est un monde de mouvement, mais où génération et corruption n'ont point de place, d'où la disparition quasi totale en physique moderne des considérations de forme\index{form} ; il n'y a pas de morphologie inanimée. Bien entendu, en biologie, il y a par contre de la morphologie. Mais alors, il n'y a plus de mathématique, au moins en tant qu'instrument de déduction. C'est cette disparition de la notion centrale de forme qui caractérise la science moderne, alors que cette notion jouait un rôle central dans l'ontologie d'Aristote."}

Ovid's\index{Ovid (Publius Ovidius Naso) (43 BCE--17 or 18 CE)} \emph{Metamorphoses} starts with the line: ``I want to speak about bodies changed into new forms." In Book I, he recounts how gods ended the status of primal chaos,\index{chaos} a ``raw confused mass, nothing but inert matter, badly combined, discordant atoms of things, confused in the one place."
In Hesiod's \emph{Theogony},\index{Hesiod (8$^{\mathrm{th}}$ c. BCE) } chaos\index{chaos} is the name of the first of the primordial deities. Then come, in that order, Gaia (Earth),  Tartarus (at the same time a place, the Deep abyss), and Eros (love).

 Chapter 2 of \emph{Structural stability and morphogenesis} is titled \emph{Form and structural stability}. It starts with the question: \emph{what is form?}\index{form} The question, in its various aspects is difficult to answer, and Thom\index{Thom, René (1923--2002)} says that it is beyond his task.
 
   Saying that a form\index{form} is a geometric figure may be a good start (if one can be satisfied with not having a definition of a ``geometric figure",\index{Euclid (c. 325--c. 270 BCE)!\emph{Elements}}\index{Aristotle (384--322 BCE)!\emph{On the soul}}\footnote{Aristotle, in \emph{On the soul}, \cite{A-Soul} 414b19, considers it useless to try to define a ``figure", and he describes the latter as a ``sort of magnitude" (425a18). We already noted that Euclid defines a figure as ``that which is contained by any boundary\index{boundary} or boundaries" (Definition 14 of Book I of the \emph{Elements}). Plato,\index{Plato (428--348 BCE)}\index{Plato (428--348 BCE)!\emph{Meno}} on the contrary, defines a figure as the boundary of a solid (\emph{Meno} 76A).} which, technically speaking, should depend on what kind of geometry we are talking about), but then, one has to introduce an equivalence relation. Topological equivalence (homeomorphism) is certainly too weak, and metric equivalence (isometry) too restrictive. For example, a homothety will also preserve form.\index{form} Thom\index{Thom, René (1923--2002)} notes that there are instances where a square drawn in a plane such that two of its sides are horizontal (and the other two vertical) does not have the same ``form" as a square placed in such a way that its sides make an angle of 45 degrees with the horizontal.  The development of the theory of forms\index{form} depends  on the use that one wants to make of it, and for Thom, the main use of this theory is in biology. Although this is the whole subject of his book \emph{Structural stability and morphogenesis}, Thom acknowledges that this is a task whose complete realization is beyond his capabilities.  
 
 Topology\index{topology} is also an adequate language for describing spaces of forms.\index{forms!space of} Thom\index{Thom, René (1923--2002)} considers equivariant Hausdorff metrics on  spaces of form. The mathematical notions of stability,\index{stability} bifurcation, attractor, singularity, universal unfolding, envelope, etc. are thoroughly used by him in this setting.

  In the article \emph{Structuralism and biology} \cite{Structuralism}, which was published the same year as the French version of   \emph{Structural stability and morphogenesis}, Thom\index{Thom, René (1923--2002)} writes that the foundations of a structure requires a precise lexicon of elementary chreods\index{chreod} and the introduction of the notion of conditional chreod, and that catastrophe theory\index{catastrophe theory} gives the mathematical models for such structures.

One of Thom's aims in his book \emph{Structural stability and morphogenesis} was to introduce in biology\index{biology} the language of differential topology,\index{topology} in particular basic notions such as differentiable manifold, vector field, genericity,\index{genericity} transversality, universal unfolding etc. In the introduction, Thom\index{Thom, René (1923--2002)} he two predecessors of him in this domain, D'Arcy Thompson\index{Thompson, D'Arcy Wentworth (1860-1948)} who we already mentioned at several occasions, and C. H. Waddington,\index{Waddington, Conrad Hal (1905-1975)}  whose concepts of ``chreod" and ``epigenetic landscape"\index{epigenetic landscape} played a germinal part in his work.  

Conrad Hal Waddington  (1905-1975)\index{Waddington, Conrad Hal (1905-1975)} , to whom Thom refers, was a well-known biologist, working on developmental biology,\index{biology} that is, the study of growth and development of living organisms.  
 The term ``chreod"\index{chreod} which he introduced in this field (from the Greek  \tg{qr'h}, which means ``it is necessary to" and  \tg{<od'os}, which means ``way") designates the transformations underwent by a cell during its development, until it finds its place as part of the organism. During this development, the cell is subject to an incredible amount of forces exerted on it from its environment to which it is permanently adjusting. An ``epigenetic landscape"\index{epigenetic landscape} is a representation of a succession of differentiation phenomena that a cell undergoes by  hills and valleys. The image represented here is extracted from Waddington's\index{Waddington, Conrad Hal (1905-1975)}  \emph{Strategy of the Genes} (1957)  \cite{Strategy} where the gene is represented as a small ball rolling in a golf field. The idea expressed by this representation is that a tiny change in the initial conditions leads to drastic changes in the way the rolling ball will take. 
   Waddington's\index{Waddington, Conrad Hal (1905-1975)}  main idea was that  the development of a cell or an embryo does not depend only on its origin, but also on the landscape that surrounds it.   In his \emph{Esquisse} (1988), Thom\index{Thom, René (1923--2002)} returns to Waddington's epigenetic landscape,\index{epigenetic landscape} describing it as a   ``metaphor which played a primordial role in catastrophe theory."\index{catastrophe theory} \cite[p.\,19]{Esquisse}    
   
       \begin{center}
\includegraphics[width=1\linewidth]{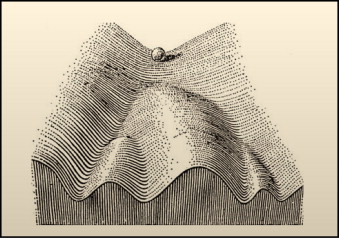} 
\vskip .1in
{\small An epigenetic landscape, from C. H. Waddington, \emph{The Strategy of the Genes}, George Allen \& Unwin, 1957, p. 29.} 
\vskip .3in
  \end{center}

  Waddington\index{Waddington, Conrad Hal (1905-1975)}  wrote two prefaces to Thom's\index{Thom, René (1923--2002)}  \emph{Structural stability and morphogenesis}, one for the first French edition, and another one for the English edition.  He writes in the first one \cite[p.\,xix]{SS}: 
  \begin{quote}\small
  I cannot claim to understand all of it; I think that only a relatively few expert topologists will be able to follow all his mathematical details, and they may find themselves less at home in some of the biology.\index{biology} I can, however, grasp sufficient of the topological concepts and logic to realise that this a very important contribution to the philosophy of science and to theoretical general biology in particular. [\ldots] Thom has tried to show, in detail and with precision, just how the global regularities with which biology deals can be envisaged as structures within a many-dimensional space. He not only has shown how such ideas as chreods,\index{chreod} the epigenetic landscape,\index{epigenetic landscape} and switching points, which previously were expressed only in the unsophisticated language of biology,\index{biology} can be formulated more adequately in terms such as vector fields, attractors, catastrophes, and the like; going much further than this, he develops many highly original ideas, both strictly mathematical ones within the field of topology,\index{topology} and applications of these to very many aspects of biology and of other sciences. [\ldots] It would be quite wrong to give the impression that Thom's book is exclusively devoted to biology. The subjects mentioned in his title, \emph{Structural stability and morphogenesis}, have a much wider reference; and he related his topological system of thought to physical and indeed to general philosophical problems. [\ldots] In biology,\index{biology} Thom\index{Thom, René (1923--2002)} not only uses topological modes of thought to provide formal definitions of concepts and a logical framework by which they can be related; he also makes a bold attempt at a direct comparison between topological structures within four-dimensional space-time, such as catastrophe hypersurfaces, and the physical structures found in developing embryos. [\ldots] As this branch of science [theoretical biology] gathers momentum, it will never in the future be able to neglect the topological approach of which Thom\index{Thom, René (1923--2002)} has been the first significant advocate.
  \end{quote}

         Another mathematician who could have been invoked in the preceding pages is Leonardo da Vinci,\index{Leonardo da Vinci (1452--1519)} who is the model---probably the supreme model---for a rare scientist-artist combination. Leonardo is also the prototype of a scholar who spent all his life learning. He was a theoretician of form.\index{form} At an advanced age, he  became thoroughly involved in  biology,\index{biology} in particular, in exploring the ideas of birth and beginning of life. He introduced some of the first known theories on the fetus. One of his notebooks is entirely dedicated to embryology.\index{embryology} He had personal ideas on the role of the umbilical cord and he developed theories on the nutritional and respiratory aspects of the embryo, as well as on its rate of change in form\index{form} during the various phases of its growth. 
         
         Leonardo\index{Leonardo da Vinci (1452--1519)} was famous for taking a long time for the execution of the works that the various patrons ordered to him, and it is interesting to know that he was blamed for spending more time on studying mathematics than on painting. Gabriel Séailles, his famous 19$^{\mathrm{th}}$ century biographer, in his book \emph{Léonard de Vinci, l'artiste et le savant : 1452-1519 : essai de biographie psychologique} \cite{Seailles}, quotes a letter from the Reverend Petrus de Nuvolaria to Isabelle d'Esté,  Duchess  of Milan, who was a leading figure of the Italian Renaissance, in which he says about Leonardo: ``His mathematical studies were, for him,  the cause of such a disgust for painting that he barely stands holding a brush."\footnote{Ses études mathématiques l'ont à ce point dégoûté de la peinture, qu'il supporte à peine de prendre une brosse.} 
Séailles also quotes Sabba da Castiglione, a writer and humanist who his contemporary,\index{Leonardo da Vinci (1452--1519)}  who writes in his memoirs: ``Instead of dedicating himself to painting, he gave himself fully to the study of geometry, architecture and anatomy."

          Leonardo\index{Leonardo da Vinci (1452--1519)} was a dedicated reader of Aristotle\index{Aristotle (384--322 BCE)}. The Renaissance was, in great part, a return to the Greek authors and in this sense, Thom\index{Thom, René (1923--2002)} was, like Leonardo, the prototype of the Renaissance man. Not only he participated in the renewed interest in Aristotle's work, but he shed a new light on it, helping us understanding better his biology\index{biology} and his mathematics.
          
          \bigskip 
          
          \noindent{\it Acknowledgements.} I am grateful to Marie-Pascale Hautefeuille who read several versions of this paper and made corrections, to Stelios Negrepontis who read an early version and shared with me his thoughts on Plato that are included in the appendix that follows this paper, and to Arkady Plotnitsky who made extremely helpful comments on an early version.  Part of this paper was written during a stay at the Yau Mathematical Sciences Centre of Tsinghua University (Beijing).

     \printindex
 
\end{document}